\documentclass{article}
\usepackage{latexsym, amsmath}
\include{amssym}
\begin{document} 
\nocite{*}
\newtheorem{thm}{Theorem} 
\newtheorem{maintheorem}{Main Theorem}   \renewcommand{\themaintheorem}{}
\newtheorem{Thm}{Theorem}
\newtheorem{prop}[thm]{Proposition}
\newtheorem{Prop}{Proposition} 
\newtheorem{lemma}[thm]{Lemma}
\newtheorem{fact}[thm]{Fact}
\newtheorem{corollary}[thm]{Corollary}
\newtheorem{cor}[thm]{Corollary}
\newtheorem{step}{Step}
\newtheorem{stp}{Step}
\newtheorem{defn}[equation]{Definition}

\newcommand\st{\mbox{$\ :\ $}}                         
\newcommand\Ind{{\rm Ind}}
\newcommand\lan{{\langle}}
\newcommand\ran{{\rangle}}
\newcommand{\order}[1]{\vert {#1} \vert}             
\newcommand\ud{{\underline{d}}}
\newcommand\ue{{\underline{e}}}
\newcommand\ux{{\underline{x}}}
\newcommand\uC{{\underline{C}}}
\newcommand\uy{{\underline{y}}}
\newcommand\eps{{\epsilon}}
\newcommand\heps{{\hat\epsilon}}
\newcommand\mustar{{\mu_*}}
\newcommand\F{{\bf F}}
\newcommand\bZ{{\bf Z}}
\newcommand\bR{{\bf R}}
\newcommand\N{{\bf N}}
\newcommand\bn{{\bf n}}
\newcommand\bs{{\bf s}}
\newcommand\cX{{\cal X}}
\newcommand\C{{\Bbb C}}
\newcommand\Q{{\Bbb Q}}
\newcommand\PD{{\bf PD}}
\newcommand\ccc{v}
\newtheorem{definition}{Definition}
    \renewcommand{\thedefinition}{} 
\newcommand\cE{{\cal E}}
\newcommand\cC{{\cal C}}
\newcommand\cEg{{\cal E}_g}
\newcommand\cM{{\cal M}}
\newcommand\cN{{\cal N}}
\newcommand\tG{{\widetilde{G}}}  
\newcommand\hG{{\widehat{G}}}  
\newcommand\hx{{\widehat{x}}}  
\newcommand\tH{{\widetilde{H}}}  
\newcommand\tV{{\widetilde{V}}}  
\newcommand\sig{\mathop{\mathrm{sig}}}
\newcommand\codim{\mathop{\mathrm{codim}}}
\newcommand\fpr{\mathop{\mathrm{fpr}}}
\newcommand\Fix{\mathop{\mathrm{Fix}}}
\newcommand{\Z}{{\Bbb Z}}
\newcommand{\PP}{\Bbb{P}^1}
\newcommand{\MMM}{{\mathcal{M}}}
\newcommand{\HHH}{{\mathcal{H}}}

\newcommand\rad{\mathop{\mathrm{rad }}}
\newcommand\Mon{\mathop{\mathrm{Mon }}}
\newcommand\Aut{\mathop{\mathrm{Aut }}}

\newcommand{\floor}[1]{\lfloor {#1} \rfloor}        
\newcommand{\ceil}[1]{\lceil {#1} \rceil}        
\newcommand\dv{{\ \mid\ }}
\newcommand\lcm{{\rm lcm}}
\newcommand\pf{\noindent {\rm Proof.\ \ \/}}
\newcommand\eop{{\hfill $\Box$ \vskip 2ex}}
\newcommand\ppp\rho
\newcommand\zetao{\zeta^\ast}
\newcommand\zzeta{\pi}

\newcommand\U{W}
\newcommand\W{W}
\newcommand\WW{W}

\title{Primitive Monodromy Groups of Genus at most Two}
\author{
Daniel Frohardt \\ 
Robert Guralnick\footnote{The second author was supported by NSF grant DMS 0653873 } \\ 
Kay Magaard}
\maketitle
\begin{abstract}
We show that if the action of a classical group $G$ on a set $\Omega$ of $1$-spaces of its natural module 
is of genus at most two, then $|\Omega| \leq 10,000$. 
\end{abstract}

\section*{Introduction}

Let $X$ be a compact, connected Riemann surface of genus $g$, and
let $\phi: X \rightarrow \PP\C$ be meromorphic of degree $n$.  Let
$B := \{ x \in \PP\C \ : \ |\phi^{-1}(x)| < n \}$ be the set of
branch points of $\phi$.  It is well known that $B$ is a finite set
and that if $b_0 \in \PP\C \setminus B$, then the fundamental group
$\pi_1(\PP\C \setminus B,b_0)$ acts transitively on $F:= \phi^{-1}(b_0)$
via path lifting.  The image of the action of $\pi_1(\PP\C \setminus
B,b_0)$ on $\phi^{-1}(b_0)$ is called the {\it monodromy group} of
$(X,\phi)$ and is denoted by $\Mon(X,\phi)$.  \\

We are interested in the structure of the monodromy group when the
genus of $X$ is less than or equal to two and $\phi$ is indecomposable
in the sense that there do not exist holomorphic functions $\phi_1
: X \rightarrow Y$ and $\phi_2 : Y \rightarrow \PP\C $ of degree
less than the degree of $\phi$ such that $\phi = \phi_2 \circ
\phi_1$. The condition that $X$ is connected implies that $\Mon(X,\phi)$
acts transitively on $F$ whereas the condition that $\phi$ is
indecomposable implies that the action of $\Mon(X,\phi)$ on $F$ is
primitive.

Our question is closely related to a conjecture made by Guralnick
and Thompson \cite{GT} in 1990. By $cf(G)$ we denote the set of
isomorphism types of the composition factors of $G$. In their paper
Guralnick and Thompson \cite{GT}  defined the set $${\mathcal E}^*(g)
= (\bigcup_{(X,\phi)} cf\Mon(X,\phi)) \  \setminus \ \{A_n, \Z/p\Z
\ : \ n > 4 \ , \ p \ \mbox{a prime} \}$$

where $X $ is a compact connected Riemann surface of genus $g$, and
$\phi : X \longrightarrow \PP(\C)$ is meromorphic, and conjectured
that ${\mathcal E}^*(g)$
is finite for all $g \in \N$. Building on work of Guralnick-Thompson
\cite{GT}, Neubauer \cite{N92}, Liebeck, Saxl \cite{LSa}, and
Liebeck, Shalev \cite{LSh}, the conjecture was established in 2001
by Frohardt and Magaard \cite{FM2}.\\

The set ${\mathcal E}^*(0)$ is distinguished in that it is contained
in ${\mathcal E}^*(g)$  for all $g$. Moreover the proof of the
Guralnick-Thompson conjecture shows that is possible to compute
${\mathcal E}^*(0)$ explicitly and indeed to describe the 
minimal covers $\phi:\PP(\C) \rightarrow |^1(\C)$ (at least those
whose monodromy group is not an alternating or symmetric of
the same degree as the cover). 

The idea of the proof of the  Guralnick-Thompson conjecture is to
employ Riemann's Existence Theorem to translate the geometric problem
to a problem in group theory as follows.  If $\phi: X \rightarrow
\PP\C$  is as above with branch points $B = \{b_1,\dots, b_r\}$,
then the set of elements $\alpha_i \in \pi_1(\PP\C \setminus B,b_0)$
each represented by a simple loop around $b_i$ forms a standardized
set of generators of $\pi_1(\PP\C \setminus B,b_0)$. We denote by $\sigma_i$ the image of $\alpha_i$ in
$\Mon(X,\phi)\subset S_F \cong S_n$. Thus we have that $$\Mon(X,\phi)
= \langle \sigma_1,\dots,\sigma_r\rangle \subset S_n $$ and that $$ \Pi_{i
= 1}^r \sigma_i = 1.$$ Moreover the conjugacy class of $\sigma_i$ in
$\Mon(X,\phi)$ is uniquely determined by $\phi$. Recall that the
index of a permutation $\sigma \in S_n$ is equal to the minimal number
of transpositions needed to express $\sigma$ as a product of such.
The Riemann-Hurwitz formula asserts that $$ 2(n+g-1) = \sum_{i=1}^r
\Ind(\sigma_i), $$ where $g$ is the genus of $X$.

\begin{defn} If $\tau_1,\dots,\tau_r \in S_n$ generate a transitive
subgroup $G$ of $S_n$ such that $\Pi_{i=1}^r \tau_i = 1$ and  $
2(n+g-1) = \sum_{i=1}^r \Ind(\tau_i) $ for some $g \in \N_0$, then
we call $(\tau_1,\dots,\tau_r)$ a {\it genus $g$-system} and $G$ a
genus $g$-group. We call a genus $g$-system $(\tau_1,\dots,\tau_r)$
primitive if the subgroup of $S_n$ it generates is primitive.
\end{defn}

If $X,\phi$ are as above, then we say that $(\sigma_1,\dots.\sigma_r)$ is
the genus $g$-system induced by $\phi$.

\begin{thm}[Riemann's existence theorem] 
For every genus $g$-system  $(\tau_1,\dots,\tau_r)$ of $S_n$ there
exists a Riemann surface $Y$ and a cover $\phi': Y \longrightarrow
\PP\C$ with branch point set $B$ such that the genus $g$-system
induced by $\phi'$ is $(\tau_1,\dots,\tau_r)$.  
\end{thm}

\begin{defn}
Two covers $Y_i,\phi_i$, $i = 1,2$ are equivalent if there exist holomorphic maps $\xi_1:Y_1\longrightarrow Y_2$ 
and 
$\xi_2:Y_2\longrightarrow Y_1$ which are inverses of one another
such that 
$\phi_1 = \xi_1\circ\phi_2$ and $\phi_2 = \xi_2\circ\phi_1$.
\end{defn}

The Artin braid group acts via automorphisms on $\Pi_1(\PP\C \setminus
B,b_0)$. We have that all sets of canonical generators of $\Pi_1(\PP\C
\setminus B,b_0)$ lie in the same braid orbit.  Also the group $G$
acts via diagonal conjugation on genus $g$-generating sets. The
diagonal and braiding actions on $g$-generating sets commute and
preserve equivalence of covers; that is if two genus $g$-generating
sets lie in the same orbit under either the braid or the diagonal
conjugation action, then the corresponding covers given by Riemann's
existence theorem are equivalent. We call two genus $g$-generating
systems {\it braid equivalent} if they are in the same orbit under
the group generated by the braid action and diagonal conjugation.
We have, see for example \cite{V} Proposition 10.14,

\begin{thm}
Two covers are equivalent if and only if the corresponding genus $g$-systems 
are braid equivalent.
\end{thm}

Suppose now that $(\tau_1,\dots,\tau_r)$ is a primitive genus
$g$-system of $S_n$. Express each $\tau_i$ as a product of a minimal
number of transpositions; i.e. $\tau_i := \prod_j \sigma_{i,j}$. The
system $(\sigma_{1,1},\dots,\sigma_{r,s})$ is a primitive genus $g$-system
generating $S_n$ consisting of precisely $2(n+g-1)$ transpositions.
By a famous result of Clebsch, see Lemma 10.15 in \cite{V}, any two
primitive genus $g$-systems of $S_n$ are braid equivalent. Thus we
see that every genus $g$- system can be obtained from one of $S_n$
which consists entirely of transpositions.

So generically we expect primitive genus $g$-systems of $S_n$ to
generate either $A_n$ or $S_n$.

We define $P{\mathcal E}^*(g)_{n,r}$ to be the braid equivalence
classes of primitive genus $g$-systems $(\tau_1,\dots,\tau_r)$
of $S_n$ such that $G:=\langle \tau_1,\dots,\tau_r\rangle$ is a
primitive subgroup of $S_n$ with $A_n \not \leq  G$. We also
define $G{\mathcal E}^*(g)_{n,r}$ to be the conjugacy classes of
primitive subgroups of $S_n$ which are generated by members of
$P{\mathcal E}^*(g)_{n,r}$.

We also define $$P{\mathcal E}^*(g) := \cup_{(n,r) \in \N^2}
P{\mathcal E}^*(g)_{n,r},$$ and similarly $$G{\mathcal E}^*(g) :=
\cup_{(n,r) \in \N^2} G{\mathcal E}^*(g)_{n,r}.$$

We note that the composition factors of elements of $G{\mathcal E}^*(g)$ are 
elements of ${\mathcal E}^*(g)$.

Our assumption that $G = \Mon(X,\phi)$ acts primitively on $F$ is
a strong one and allows us to organize our analysis along the lines
of the Aschbacher-O'Nan Scott theorem exactly as was done in the
original paper of Guralnick and Thompson \cite{GT}. We recall the
statement of the Aschbacher-O'Nan-Scott theorem from \cite{GT}

\begin{thm} 
Suppose $G$ is a finite group and $H$ is a maximal subgroup of $G$
such that $$\bigcap_{g \in G} H^g = 1.$$ Let $Q$ be a minimal normal
subgroup of $G$, let $L$ be a minimal normal subgroup of $Q$, and
let $\Delta = \{L=L_1, L_2, \dots, L_t\} $ be the set of $G$-conjugates
of $L$. Then $G= HQ$ and precisely one of the following holds:

\begin{itemize}
\item[(A)] $L$ is of prime order $p$.
\item[(B)] $F^*(G) = Q \times R$ where $Q \cong R$ and $H \cap Q = 1$.
\item[(C1)] $F^*(G) = Q$ is nonabelian, $H \cap Q = 1$.
\item[(C2)] $F^*(G) = Q$ is nonabelian, $H \cap Q \neq 1 = L \cap H$.
\item[(C3)] $F^*(G) = Q$ is nonabelian, $H \cap Q = H_1 \times \dots
\times H_t$, \\ where $H_i = H \cap L_i \neq 1$, $1 \leq i \leq t.$
\end{itemize} 
\end{thm}

We summarize briefly what is known about $G{\mathcal E}^*(0)$ and
$P{\mathcal E}^*(0)$.  The members of $G{\mathcal E}^*(0)$ that
arise in case (C2) were determined by Aschbacher \cite{A}. In all
such examples $Q = A_5 \times A_5$.  Shih \cite{Shi} showed that
no elements of $G{\mathcal E}^*(0)$ arise in case (B) and Guralnick
and Thompson \cite{GT} showed the same in case (C1). In his thesis
Neubauer \cite{N} showed that in case (A) either $G''=1$ and $G/G'$
is an abelian subgroup of $GL_2(p)$, or that $n \leq 256$.  Recently
Magaard, Shpectorov and Wang \cite{MSW}, determined all elements
of $P{\mathcal E}^*(0)_{n,r}$ with $n \leq 256$. The elements $G$
of $G{\mathcal E}^*(0)$ arising in case (C3)  have generalized
Fitting subgroups with fewer than $5$ components; i.e. $t \leq 5$.
This was shown by Guralnick and Neubauer \cite{GN} and later
strengthened by Guralnick \cite{G} to $t \leq 4$. Moreover Guralnick
showed that the action of $L_i$ on the cosets of $H_i$ is a member
of $G{\mathcal E}^*(0)$. In case (C3) where $L_i$ is of Lie type
of rank one all elements of $G{\mathcal E}^*(0)$ and $G{\mathcal
E}^*(1)$ were determined by  Frohardt, Guralnick, and Magaard
\cite{FGM1}, moreover they show that $t \leq 2$. In \cite{Kong} Kong shows that
if $G$ is an almost simple group of type $L_3(q)$, 
then $G \in G{\mathcal E}_{(q^2+q+1,r)}^*(g)$ with
$g \leq 2$ only if $q \leq 13$, and $G \in G{\mathcal E}_{(q^2+q+1,r)}^*(0)$ if and only if
$q \leq 7$. 
Combining the results of Frohardt, Magaard
\cite{FM3} with those of Liebeck, Seitz \cite{LSe} we have that if
$F^*(G)$ is exceptional of Lie type and $G \in G{\mathcal E}^*(0)_{n,r}$,
then $n \leq 65$. In \cite{GS} Guralnick and Shareshian show that
$G{\mathcal E}^*(0)_{n,r} = \empty$ if $r \geq 9$.  Moreover they
show that if $G \in G{\mathcal E}^*(0)_{n,r}$ with $F^*(G)$ alternating
of degree $d < n$, then either $r \leq 4$ or $r = 5$ and $n =
d(d-1)/2$.
In \cite{M} Magaard
showed that if $F^*(G)$ is sporadic and $G \in G{\mathcal E}^*(0)_{n,r}$
then  $n \leq 280$. We would like to take this opportunity to point out that 
$\Aut(HS) \in G{\mathcal E}^*(0)_{100,3}$ which was missed in cite \cite{M}. Furthermore we
thank the referee for pointing out that $\Aut(HS)$ possesses four genus zero systems in its action 
on $100$ points with signatures $(2, 4, 10), (2, 5, 6), (2, 4, 5)$, and $(2, 4, 6).$
The referee has further pointed out that first two of these genus zero systems are  
rational and rigid. This is because in both of these cases the involution has an odd number 
of transpositions, and therefore the corresponding genus 0 field is rational. 
Thus there exists $\phi: \PP\Q\rightarrow \PP\Q$ of a degree $100$ with
monodromy group $\Aut(HS)$.

This leaves open the cases $F^*(G) = A_d^t$, $t,r \leq 4$ and the
cases $F^*(G) = L^t$, $t \leq 4$ where $L$ is a classical group of Lie type. 
In light of the results of \cite{AGM} we suspect that if $G$ is in the second
case and $G \in G{\mathcal E}^*(0)$, then $L_i/H_i$ is a point
action, i.e. equivalent to an action of $L_i$ acting on an orbit of
one-spaces of its natural module. Hence they are the focus of this
paper.

Another problem closely related to the Guralnick-Thompson conjecture is
the  description of the monodromy groups from the generic Riemann surface  of genus
$g$ to $\PP(\C)$ of degree $n$.   This is related to Zariski's thesis where
he answered a conjecture of Enrique by showing that the generic
Riemann surface of genus $g > 6$ does not admit a solvable map of fixed degree $n$ to
$\PP(\C)$ (i.e. where the monodromy group is solvable).  The condition on 
$n$ being fixed was removed in \cite{GN}.  Note that any Riemann surface of genus
at $6$ admits a degree $4$ map to $\PP(\C)$ (and so is solvable).   Interestingly,
Zariski's methods were mostly group theoretic. 

  Recall that the images
of the canonical generators of $\pi_1(\PP\C \setminus B,b_0)$ are
determined uniquely up to conjugacy in $G$. We say that a $G$-cover
of $\PP\C$ has \emph{ramification type} $C_1,\dots, C_r$ if the
$i$'th canonical generator lies in conjugacy class $C_i$ of $G$.
The moduli space of $G$-covers of $\PP\C$ with ramification type
$C_1,\dots, C_r$ is a \emph{Hurwitz space} and is denoted by
$\HHH(G,0,C_1,\dots,C_r)$. Via the Riemann-Hurwitz formula we see
that every $G$-cover $X \in \HHH(G,0,C_1,\dots,C_r)$ has the same
genus $g$. So the forgetful functor ${\mathcal F} :
\HHH(G,0,C_1,\dots,C_r) \rightarrow \MMM_g$ is well defined and so
the problem of determining maps of degree $n$ from the generic
Riemann surface of genus $g$  can be rephrased as follows:

For which groups $G$ and which ramification types $C_1, \dots, C_r$
of $G$ is the forgetful functor ${\mathcal F} : \HHH(G,0,C_1,\dots,C_r)
\rightarrow \MMM_g$ dominant; i.e. is the image of $\HHH(G,0,C_1,...,C_r)$
dense in $\MMM_g$?

Now Theorem 2 of Guralnick Magaard \cite{GM} shows that if the image
of $\HHH(G,0,C_1,...,C_r)$ under the forgetful functor is dense in
$\MMM_g$, then one of the following holds 
\begin{enumerate}
\item $g \leq 2$,
\item $g = 3$ and $G$ is affine of degree $8$ or $16$, 
\item $g = 3$ and $G \cong L_3(2)$, \item $g \geq 3$ and $G \cong S_n, n \ge (g+2)/2$ or
$A_n, n > 2g$.  
\end{enumerate}

It is well known that $S_n$ does cover $\PP\C$ generically. However
it was only in 2006 when Magaard and V\"olklein \cite{MV} proved
that $A_n$ and $L_3(2)$ also cover $\PP\C$ generically.  It was
later shown by Magaard, V\"olklein and Wiesend \cite{MVW} that
$AGL_3(2)$ and $AGL_4(2)$ cover $\PP\C$ generically. This leaves
only the first possibility, and is a reason why our ultimate goal
is to determine $P{\mathcal E}^*(g)$ where $g \leq 2$.

Our two primary results here are Theorem \ref{main.result}, which shows that if $n > 10^4$
then the elements of $P{\mathcal E}_{n,r}^*(g)$ with $g \leq 2$
are not point actions of classical groups, and Theorem \ref{basic.result} which is 
more technical but can be applied to a wider class of actions.
Combining Theorem \ref{main.result} with the main theorem of \cite{FGM2} shows that if $n > 10^4$,
then the elements of $P{\mathcal E}_{n,r}^*(g)$ with $g\leq 2$ are generally not subspace actions of 
classical groups. The potential exceptions to the statement 
are also explicitly given in \cite{FGM2}. These potentially exceptional actions are precisely those actions whose 
permutation modules do not contain the permutation module of
the action on singular points as a submodule. 
The main result of \cite{AGM} determines all classes of maximal subgroups of classical groups whose permutation
module does not contain the permutation module of the action on singular points. For these classes of 
maximal subgroups we hope to establish the hypotheses of Theorem \ref{main.result} which would then show that 
if $n > 10^4$, then the elements of $G{\mathcal E}_{n,r}^*(g)$ with $g \leq 2$ are either cyclic of prime
order $n$ or contain the alternating group $A_n$.

To establish Theorem \ref{basic.result} we show that for any pair $(G,\Omega)$,
where $G$ is a classical group acting primitively on a set $\Omega$ such that the hypotheses of 
Theorem \ref{basic.result} are satisfied, and 
any generating $r$-tuple $(\tau_1, \dots, \tau_r)$ of $G$ which satisfies the product $1$ condition, then the expression 
$\sum_{i=1}^r \Ind(\tau_i)$ is greater than $(2+\epsilon)n$ for some positive constant $\epsilon.$
We achieve this by proving effective lower bounds on $\Ind(\tau_i)$ using Scott's Theorem \ref{generation.fact}
and the technique of translating tuples, see Lemma \ref{magnus}.

\section{Statement of Results}

\begin{definition}
$\ux = (x_1, x_2, \ldots, x_r)$ is a {\em normalized generating $r$-tuple\/}
for $G$ provided
\begin{enumerate}
\item $G = \lan x_1, x_2, \ldots, x_r \ran$
\item $x_1 x_2 \ldots x_r = 1$
\item $x_i \neq 1$, $i = 1, 2, \ldots, r$
\end{enumerate}

If, in addition, 
$G$ is a transitive permutation group of degree $N$ and
$$\sum \Ind(x_i) = 2(N+g -1)$$
then $\ux$ has genus $g$.
\end{definition}

The formula above is the Riemann-Hurwitz Formula (RH).  
The Riemann Existence Theorem \cite{GT} 
guarantees that given a normalized generating tuple $\ux$ 
for a permutation group $G$ there is a surface $X$ and a 
covering $\rho : X \mapsto \PP(\C)$ such that $G \cong \Mon(X,\rho)$ 
and the genus of $X$ is the genus of the tuple $\ux$, written $g(\ux)$.

The primary result of this paper is the following.

\begin{thm}
\label{main.result}
If $(G,\Omega)$ is a primitive classical point action of degree 
at least $10^4$, then the action has genus larger than 2.
\end{thm}

The case of point actions will lead almost all the examples
(indeed using \cite{FGM} and some ongoing work of {AGM},
one can eliminate most other situations). 

\bigskip
The proof of Theorem~\ref{main.result}
uses inequalities based on RH and estimates for the fixed point 
ratios of elements of $G$.

\begin{definition}
For $x$ a permutation of the finite set $\Omega$, let
$F_\Omega(x)$ (or $F(x)$) denote the fixed points of $x$ on
$\Omega$ and let $f_\Omega(x)$ (or $f(x)$) denote the
{\em fixed point ratio\/} of this permutation.
That is, $f(x) = \order{F(x)}/\order{\Omega}$.
\end{definition}

\begin{definition}
Let $V$ be a vector space and let $x \in \Gamma L(V)$.
If $x$ acts as a permutation on the set $\Omega$ then 
the triple $(x,V,\Omega)$
satisfies Grassmann Condition $\eps$ provided 
$$f_\Omega(x) < \frac{\order{\W}}{\order{V}} + \eps$$
for some eigenspace $\W$ for the action of $x$ on $V$.

A classical group $G$ with natural module $V$ acting as a permutation
group on the set $\Omega$ satisfies Grassmann Condition $\eps$
provided $(x,V,\Omega)$ sastifies Grassmann Condition $\eps$
for every $x \in G$.
\end{definition}

Note:  For the purposes of the previous definition, an eigenspace 
for the action of $x$ on $V$ is a set $\W \subset V$ which is a subspace 
of $V$ over some (possibly proper) subfield of the field of definition 
on which $x$ acts as a scalar.  Note that $|\W|$ does not depend on 
its field of scalars.

\bigskip

The role of Grassmann Conditions in the proof of Theorem~\ref{main.result} 
is apparent in the statement of the following technical results which 
together yield Theorem~\ref{main.result}.
The key feature of the point actions is that with known exceptions 
they satisfy Grassmann Condition $1/100$.  Theorem~\ref{basic.result} 
also applies to other actions that satisfy this condition.

\begin{thm}
\label{basic.result}
Let $G$ be a linear group with module $V$ where $V$ contains
at least $10^4$ projective points and no constituent for the action 
of $G$ on $V$ has dimension $1$.  
If $\ux$ is a normalized generating $r$-tuple for $G$ in some primitive
permutation action, then one of the following is true.
\begin{enumerate}
\item $g(\ux) > 2$.
\item 
$G$ does not satisfy Grassmann Condition $1/100$.  
More specifically, for some $i \in \{1, \ldots, r\}$, 
the group $\langle x_i \rangle$ contains an element $y$ 
that violates Grassmann Condition $1/100$.  
\item The characteristic of $V$, the dimension of $V$ over its prime field,
and the signature of $\ux$ are given in Table~\ref{pnd.table}.
\end{enumerate}
\end{thm}

Note:
The {\em signature\/} $\sig(\ux)$ of the $r$-tuple 
$\ux = (x_1, x_2, \ldots, x_r)$
is the $r$-tuple $(d_1,d_2, \ldots, d_r)$,
where $d_i = o(x_i)$.


\begin{thm}
\label{grassmann.exceptions}
Let $G$ be a classical group with natural module $V$.
Let $\Omega$ be a primitive point action for $G$ 
with $\order{\Omega} \geq 10^4$ and
assume that $G$ does not satisfy Grassmann Condition $1/100$.
Then $g(\ux) > 2$ for every normalized generating tuple $\ux$ such 
that $\langle x_i \rangle$ contains an element $y$ violating 
Grassmann Condition $1/100$.
\end{thm}

\begin{thm}
\label{touch.up}
Let $G$ be a classical group with natural module $V$ 
Assume $\ux$ is a normalized generating tuple for $G$ and that $\Omega$ 
is a primitive point action for $G$
with $\order{\Omega} \geq 10^4$.
If the characteristic of $V$, the dimension of $V$ over its prime field,
and the signature of $\ux$ are given in Table~\ref{pnd.table} then
$g(\ux) > 2$.
\end{thm}



\begin{table}
\caption{Characteristic, Dimension and Signature 
of Exceptional Cases in Theorem~\ref{basic.result}}
\label{pnd.table}
\begin{center}
\begin{tabular}{ccc}
$p$ & $\dim_{\F_p}(V)$ & $\sig(\ux)$ \\
\hline
$11$ & $5,6$ & $(2,3,7)$\\
$7$ & $6$ & $(2,3,7)$\\
$5$ & $7,8,9$ & $(2,3,7)$\\
$3$ & $12$ & $(2,3,7)$\\
$2$ & $14,15, \ldots, 21$ & $(2,3,7)$\\
$3$ & $10$ & $(2,3,8)$\\
$2$ & $16$ & $(2,4,5)$ 
\end{tabular}
\end{center}
\end{table}

\begin{definition}
The almost simple classical group $G$ has a {\em point action\/} on $\Omega$
provided $G$ has a natural module $V$
of dimension $n$ over $\F_q$ where $(G,\Omega,n,V)$ satisfy one of the following conditions.

\begin{enumerate}
\item[$L$]: 
$F^\ast(G) \cong L_n(q)$, and $\Omega$ is the set of all
points in $V$.  $n \geq 2$. 

\item[$O^\eps,\bs$]:
$F^\ast(G) \cong O^\eps_{n}(q)$, $V$ is a non-degenerate orthogonal space of type
$\eps$, and $\Omega$ is the set of singular points in $V$.
$n$ is even, $n \geq 6$, $\eps = +$ or $-$.

\item[$O^\eps,\bn$]:
$F^\ast(G) \cong O^\eps_{n}(q)$, $V$ is a non-degenerate orthogonal space of type
$\eps$, and $\Omega$ is the set of $+$-type points in $V$.
$n$ is even, $n \geq 6$, $\eps = +$ or $-$.

\item[$O,\bs$]:
$F^\ast(G) \cong O_n(q)$, $V$ is a non-degenerate orthogonal space, 
and $\Omega$ is the set of singular points in $V$.
$n$ is odd, $n \geq 5$, and $q$ is odd.

\item[$O,\delta$]:
$F^\ast(G) \cong O_n(q)$, $V$ is a non-degenerate orthogonal space, 
and $\Omega$ is the set of $\delta$-type points in $V$.
$n$ is odd, $n \geq 5$, $\delta = +$ or $-$, and $q$ is odd.

\item[$Sp$]:
$F^\ast(G) \cong Sp_{n}(q)$, $V$ is a non-degenerate symplectic space
and $\Omega$ is the set of points in $V$.
$n$ is even, $n \geq 4$.

\item[$Sp,\delta$]:
$F^\ast(G) \cong Sp_{n}(q)$, 
$V$ is a symplectic space, and $\tV$ is an orthogonal space of
dimension $n+1$ such that $\rad \tV$ is anisotropic
of dimension $1$ and $V \cong \tV/\rad \tV$, 
and $\Omega$ is the set of all complements to 
$\rad \tV$ in $\tV$ of type $\delta$.
$n$ is even, $n \geq 4$, $\delta = +$ or $-$, and $q$ is even.

\item[$U,\bs$]:
$F^\ast(G) \cong U_n(q^{1/2})$, $V$ is a non-degenerate hermitian space,
and $\Omega$ is the set of singular points in $V$.
$n \geq 3$, $q$ is a square,

\item[$U,\bn$]:
$F^\ast(G) \cong U_n(q^{1/2})$, $V$ is a non-degenerate hermitian space,
and $\Omega$ is the set of nonsingular points in $V$.
$n \geq 3$, $q$ is a square,
\end{enumerate}
\end{definition}

\smallskip
We prove Theorems~\ref{basic.result}, \ref{grassmann.exceptions}, and \ref{touch.up} 
in the subsequent sections.
Since the action of a classical group $G$ on its natural module $V$ 
satisfies the hypotheses of Theorem~\ref{basic.result},
it is evident that Theorem~\ref{main.result} follows from these theorems.

%
%

\section{Proof of Theorem~\ref{basic.result}}

\subsection{Notation and preliminary results}
\label{notation.section}

Let $G$ be an almost simple classical group with natural
module $V$ of dimension $n_q$ over $\F_q$ and let $p$ be 
the characteristic of $\F_q$. 
Then $V_{\F_p}$ 
is an $\F_p$-vector space and all elements of $G$
correspond to $\F_p$-linear maps.  
We have $G = \hG/Z$ where $\hG \subseteq GL(V_{\F_p})$ and
$Z$ acts as scalars on $V_{\F_p}$.
Set $n_p = \dim_{\F_p} V_{\F_p}$, so that $n_p = n_q \log_p(q) $.

\begin{definition} 
$v_q(y)$ [resp., $v_p(y)$] is the codimension of the largest eigenspace
of the action of an associate of $y$ on $V$ [resp., $V_{\F_p}$].
\end{definition}
Regarding $V$ as an $\F_p$-space, $v_p(x) = \max(\codim C_V(\hx) 
\ : \ \hx \mapsto x$ under $\hG \rightarrow G)$.

Let $\Omega$
be a primitive $G$-set of order $N$.
Let $\ux = (x_1,\ldots, x_r)$ be a normalized generating tuple for $G$.

Let $g = g(\ux)$, and let $$\ud = (d_1, \ldots, d_r)$$ be the signature of $\ux$, so
that
$d_i = o(x_i), i = 1, \ldots, r$.


When the context is clear, we will write $n$ instead of $n_q$ or $n_p$ and
$v$ instead of $v_q$ or $v_p$.


\bigskip
%
The Cauchy-Frobenius Formula says that
if $x \in G$ has order $d$, then 
$$\Ind(x) = N - \frac1d \sum_{y \in \lan x \ran} F(y).$$ 

Combining this with (RH), we have
\begin{equation} \label{CF.consequence}
\sum_{i=1}^r \frac1{d_i} 
    \left( 1 + \sum_{y \in \lan x_i \ran^\sharp} f(y) \right)
    = r - 2 - 2\left(\frac{g-1}N\right).
\end{equation}

\begin{definition}
$$\eps_0 = 2\left(\frac{g-1}N\right),$$ 
$$A(\ud) = \sum\frac{d_i-1}{d_i}.$$
\end{definition}

\begin{definition}
For $x \in G$, with $o(x) = d$, set 
$$\kappa(x) = \frac1d\left(1+ \sum_{y \in \lan x \ran^\sharp} p^{-v(y)}\right)$$
\end{definition}

\begin{fact}
\label{kappa.inequality}
If $G$ satisfies Grassmann condition $\eps$ 
then 
$$\sum \kappa(x_i) > r-2 - A(\ud)\eps - \eps_0$$
\end{fact}

\pf 
Since $p^{-v(y)} = \frac{|W|}{|V|}$ where $W$ is the largest 
eigenspace for $V$, 
we have $f(y) < p^{-v(y)} + \eps$ for all $y \in G$.
Therefore, 
$$r - 2 - \eps_0 < \displaystyle \sum_{i=1}^r \frac1{d_i}
\left( 1 + (d_i-1) \eps - \sum_{y \in \lan x \ran^\sharp} p^{-v(y)}\right)
< A(\ud) \eps + \sum \kappa(x_i).$$
\eop

The relevance of this result can be seen from the main result of \cite{FM1}.
\begin{thm}[Grassmann Theorem]
There is a function $\heps : \N \rightarrow \bR^+$
such that 
   \begin{enumerate}
      \item $(G,\Omega)$ satisfies Grassmann condition $\heps(m)$ 
         whenever $(G,\Omega)$ is a classical subspace action of degree $m$,
         and
      \item
         $\displaystyle \lim_{m\rightarrow \infty} \heps(m) = 0$.
   \end{enumerate}
\end{thm}


In the balance of this subsection we obtain upper bounds for $\kappa(x)$
that will be used in the proof of Theorem~\ref{basic.result}. 

Set
$$\zeta(d) = \zeta(d,p) = \frac1d\left( 1 + \sum_{m|d,m>1} \phi(m) p^{-1}\right),$$
where $\phi$ is the Euler $\phi$-function on integers. When $a$ is not an integer
we take $\phi(a) = 0$.
Since
\begin{eqnarray*}
\kappa(x) & = &\frac1d\left(1 + \sum_{m|d,m > 1} \phi(m) p^{-v(x^{d/m})}
\right) \\
& = &\frac1d\left(1 + \sum_{m|d, m<d} \phi(\frac{d}{m})p^{-v(x^m)}\right),
\end{eqnarray*}
it follows that
if $x$ has order $d$, then 
\begin{equation} \kappa(x) \leq \zeta(d) = \frac1d + \frac1p - \frac1{dp}.
\end{equation}
Note that $\zeta$ is a decreasing function of both $d$ and $p$.

For each positive integer $s \geq 1$, set
$$\zeta_s(d) = \frac1d\left( 1 + \phi(d)\cdot p^{-s}
+ \sum_{m|d, 1 < m < d} \phi(m) p^{-1}\right).$$
More generally, for a finite sequence $s_1,s_2,\ldots, s_l$ of positive
integers, let
$$\zeta_{s_1,s_2,\ldots,s_l}(d) = \frac1d\left( 1 
+ \sum_{i=1}^l \phi(d/i) p^{-s_i}
+ \sum_{m|d,1 < m < d/l} \phi(m) p^{-1}\right).$$

The following statement is evident.
\begin{fact}
\label{zeta.sub.s}
If $x$ has order $d$ and $v(x^{i}) \geq s_i, i = 1, \ldots, l$, then
$\kappa(x) \leq \zeta_{s_1,\ldots,s_l}(d)$.
\end{fact}

The estimates for $\kappa(x)$ can be further refined by taking into consideration
the possible actions of elements of a given order on a vector space over $\F_p$.

\begin{definition}
For each prime $p$ and integer $d \geq 2$ let $\mustar(d,p)$
be the smallest positive integer $\mu$ such that 
$\mu = \dim{([V,x])}$
for some linear operator $x$
of order $d$ acting on a vector space $V$ over $\F_p$.
\end{definition}

Note that each $x \in G$ is the image of some element $\hx$ in $\hG$ 
with $\dim C_{V_p}(\hx) = v(x)$.

If $y \in G$ has order $m$ then $v(y) \geq \mustar(m,p)$.
This inequality holds in particular when
$m | d$, $o(x) = d$, and $y = x^{d/m}$. 
Set
$$
\zetao(d) = \zetao(d,p)  = 
          \frac1d \left( 1 + \sum_{m | d, m > 1} \phi(m) p^{-\mustar(m,p)}\right).  
$$
Then
\begin{equation}
\kappa(x) \leq \zetao(d).
\end{equation} 

Similarly, if 
$$
\zetao_{s_1,\ldots,s_l}(d) = 
          \frac1d \left( 1 + \sum_{m | d, m > 1} \phi(m) p^{-\alpha(d/m)}\right),
\ \  \alpha(i) = \max(s_i,\mustar(d/i,p)),
$$
then
\begin{equation} 
\kappa(x) \leq \zetao_{s_1,\ldots, s_l}(d)
\end{equation} 
whenever
$v(x^i) \geq s_i$, $i =1, \ldots, l$.

\begin{lemma}
\label{Lemma5}
\begin{enumerate}
\item
If $p > 2$ then $\zetao(d) < \frac3d + .04$.
\item
If $p = 2$ then $\zetao(d) < \frac4d + .032$.
\end{enumerate}
\end{lemma}

\pf
Suppose $p > 3$.
Then $\mustar(d) = 1$ if and only if $d=p$ or $d | p-1$,
and $\mustar(d) > 1$ for all other $d$.  Since at most
$p-1$ nontrivial powers of an element have order $p$
and at most $p-2$ nontrivial powers of an element have
order  dividing $p-1$
this implies that 
$\zetao(d,p) \leq \frac1d(1 + (2p-3)p^{-1} +  (d-1-(2p-3))p^{-2})
< \frac1d(1 + 2 + d/p^2) =  3/d + 1/p^2 \leq 3/d + 1/5^2$.
If $p = 3$, then $\mustar(m,p) = 1$ if and only if $m = 2$ or $3$,
and  $\mustar(m,p) = 2$ if and only if $m = 4,6$, or $8$.
This implies that
$\displaystyle \sum_{m|d,\mustar(m,p) = 2} \phi(m) \leq 
\phi(4) + \phi(6) + \phi(8) = 8$.
Therefore $\zetao(d,3) \leq \frac1d(1 + 3 \cdot 3^{-1} + 8 \cdot3^{-2}
+ (d-12) \cdot 3^{-3}) < 3/d + 1/27$.

For $p=2$, we note that 
$\mustar(m,2) = 1$ if and only if $m=2$;
$\mustar(m,2) = 2$ if and only if $m=3$ or $4$;
$\mustar(m,2) = 3$ if and only if $m=6$ or $7$; and
$\mustar(m,2) = 4$ if and only if $m=5,8,12,14$, or $15$.
It follows from this that $\zetao(d,2) \leq 4/d + 1/32$.
\eop

\begin{corollary}
\label{zeta.lemma}
Let $x \in G$ have order $d$, and let $k$ be a real number.
\begin{enumerate}
\item
\label{3d.bound}
If $p > 2$
and $\zeta(d) \geq k > .04$ then $d \leq \displaystyle \frac3{k-.04}$.
\item
If $p = 2$ and $\zeta(d) \geq k > .032$ then $d \leq \displaystyle \frac4{k-.032}$.
\end{enumerate}
\end{corollary}

The precise value of $\mustar(d,p)$, the smallest possible commutator dimension for
an element of order $d$ over $\F_p$, can be computed using the following statement.

\begin{fact}
\label{mu.fact}
\begin{enumerate}
\item  \label{decomp.statement}
If $d_p$ is the largest power of $p$ dividing $d$ and
$d_{p'} = d/d_p$, then
$\mustar(d,p) = \mustar(d_p,p) + \mustar(d_{p'},p)$.
\item \label{jordan.statement}
For $a \geq 1$, $\mustar(p^a,p) = p^{a-1}$.
\item
If $(d,p) = 1$ then either
$\mustar(d,p)$ is the exponent of $p \pmod d$ or
$\mustar(d,p) = \mustar(a,p) + \mustar(b,p)$
for some integers $a,b$ with $ab = d$, $a, b > 1$,
and $(a,b) = 1$.
\end{enumerate}
\end{fact}

\pf
We may assume that $d > 1$.
Suppose $x$ is an operator of order $d$ on $V$ that achieves
the minimum commutator dimension.
Without loss, assume that $\dim V$ is minimal.

$V$ is a direct sum of indecomposable $\F_p\lan x\ran$-submodules
$V_i$.  Setting $x_i = x\!\!\mid_{V_i}$ we have
$o(x) = \gcd(\{o(x_i)\})$ and $\dim ([V,x]) = \sum \dim([V_i,x_i])$.
Since $\dim ([V_i,x_i^{m}]) \leq \dim ([V_i,x_i])$ for all
$m \in \N$, by  minimality of $\dim ([V,x])$ we may
assume that $o(x_i)$ is relatively prime to $o(x_j)$
when $i \neq j$.

To prove statement~\ref{jordan.statement}, suppose $d = p^a$.  
Then $V$ consists
of a single Jordan block with eigenvalue $1$.
The order of a Jordan block of size $b$ with 
eigenvalue $1$ is $p^a$ where $p^{a-1} < b \leq p^a$.
[Proof:  $(y-1)^{b-1} \neq 0$ and $(y-1)^b = 0$ 
imply $y^{p^k} = 1$ exactly when $p^k \geq b$.]
Therefore $p^a \geq \dim V > p^{a-1}$, whence $\dim V = p^{a-1} + 1$
by minimality, and $\mustar(a,p) = \dim ([V,x]) = \dim V - 1 = p^{a-1}$.

To prove \ref{decomp.statement}, note that since $ab \geq a-1+b$
for positive integers $a$ and $b$, for unipotent $u$ and semisimple
$s$ the commutator dimension of $u \otimes s$ is always at
least as large as the commutator dimension of $u \oplus s$.

The last statement follows easily from the fact that if
$x$ acts irreducibly and semisimply on $V$ then $\dim V$
is the exponent of $p \pmod d$.
This completes the proof of ({\bf\ref{mu.fact}}).
\eop

\subsection{System bounds}
\label{scottbound.section}

The results of the previous subsection apply to individual elements.  
We shall require stronger bounds, which depend on the system, not merely the 
individual generating elements.  
As in \cite{FM2},
we use a result of L.~Scott 
on linear groups together with a fact
about group generation to control the contributions of elements
with large fixed point ratios to the index sum.

\begin{thm} [Scott]
\label{generation.fact}
Suppose $\hG$ is a group of linear operators on $V$ with $[V,\hG] = V$ 
and $C_V(\hG) = 1$.  
If $\hG = \lan g_1, \dots, g_r \ran$ where  $\prod g_i = 1$, then
$\sum \dim ([V,g_i]) \geq 2 \dim V$.
\end{thm}

\pf
See \cite{Scott}.
\eop

\begin{lemma} 
\label{magnus}
Assume that $\ue$ is an ordered $r$-tuple that is a permutation
of one of the following.
\begin{enumerate}
\item $(m,m,1, \ldots, 1)$, $m \geq 1$.
\item $(2,2,m,1, \ldots, 1)$, $m \geq 2$.
\item $(2,3,m,1, \ldots, 1)$, $m = 3, 4$, or $5$.
\end{enumerate}
Set $\displaystyle C_i = C_i(\ue) = \frac2{e_i ( 2 - A(\ue) ) }$.

Let $H$ be a group with generators $\{y_1, \ldots, y_r \}$ 
where $y_1y_2\cdots y_r =1$. 
Then there is an ordered $M$-tuple $(z_1, \ldots, z_M)$,
of elements 
of $H$, where $M = \sum_j C_j$ such that the following conditions hold.
\begin{enumerate}
\item
$z_1 z_2 \cdots z_M = 1$
\item
$\displaystyle \{1, \ldots, M \} = 
\cup_{i=1}^n {\cal C}_i$ (disjoint union) where
$|{\cal C}_i| = C_i$ and $z_j$ is conjugate to $y_i^{e_i}$ for 
all $j \in {\cal C}_i$.
\item
The group $K$ generated by $\{ z_j \}$ is normal of index $2/(2-A(\ue))$ 
in $H$, and $H/K$ is cyclic, dihedral, or isomorphic to $Alt_4$, $Sym_4$, or
$Alt_5$. 
\end{enumerate}
\end{lemma}

\pf
By well-known properties of generators and relations (see
\cite{Magnus}, for example), if $\ue$ is one of the specified tuples,
then the group $\lan y_i, i = 1, \ldots, r \ | \ y^{e_i} = \prod_i y_i = 1
\ran$ is, in the respective cases,
cyclic of order $m$, dihedral of order $2m$, or isomorphic to $Alt_4$, $Sym_4$, or
$Alt_5$. 
In each case, this group has order $2/(2-A(\ue))$.  
The statements follow from the proof of Lemma~3.2 in \cite{FM1} or 
from \cite{GN}.
\eop

Note that if $\uC(\ue) = (C_1(\ue),\dots,C_r(\ue))$ then
$$\uC(m,m,1, \ldots, 1) = (1,1,m, \ldots, m)$$
$$\uC(2,2,m,1, \ldots, 1) = (m,m,2,2m \ldots, 2m)$$
$$\uC(2,3,3,1, \ldots, 1) = (6,4,4,12, \ldots, 12)$$
$$\uC(2,3,4,1, \ldots, 1) = (12,8,6,24, \ldots, 24)$$
$$\uC(2,3,5,1, \ldots, 1) = (30,20,12,60, \ldots, 60)$$

Assume now that $\ux$ is a normalized generating $r$-tuple for $G$,
a classical group with natural module $V$ with $\dim (V/C_\hG(V)) = n$.

\begin{lemma}
\label{translation.lemma}
If $\ue$ and $\uC$ are as above,
then, for each $i^\ast$ in $\{ 1, \dots, r\}$, 
$$(C_{i^\ast}-1)v(x_{i^\ast}^{e_{i^\ast}}) 
    + \sum_{i\neq {i^\ast}} C_i v(x_i^{e_i}) \geq n.$$
If $p = 2$, then 
$$\sum C_i v(x_i^{e_i}) \geq 2n.$$
\end{lemma}

\pf
%
We apply Theorem~\ref{generation.fact} to the preimages ${\hat{z}_j}$ 
of the elements $z_j$ under the homomorphism $\hG \rightarrow G$.
In general, we can choose $M-1$ preimages $\hat{z}_j$ 
so that 
$\dim ([V,\hat{z}_j]) = v(x_i^{e_i})$, when $j \in \cC_i$.
If $j^\ast$ is the remaining subscript and $j^\ast \in \cC_{i^\ast}$, then
$\dim([V,\hat{z}_{i^\ast}]) \leq n$, and we have the first 
statement.  

If $p= 2$, then 
$\dim([V,\hat{z}_j]) = v(x_i^{e_i})$ 
whenever $j \in \cC_i$ because $|{\bf F}^\times| = 1$.  
\eop

\begin{fact}
\label{zetas}
Suppose $r = 3$ and $d_1 \leq d_2 \leq d_3$.
\begin{enumerate}
   \item If $n > d_1$, then $v(x_i) \geq 2$ for $i \geq 2$.
   \item If $n > d_2$, then $v(x_1) \geq 2$ for all $i$.
   \item If $n \geq 4$ and $d_1 \leq 3$, then $\kappa(x_i) < \zeta_2(d_i)$ for $i > 1$.
   \item If $n \geq 4$ and $d_2 \leq 3$, then $\kappa(x_i) < \zeta_2(d_i)$ for all $i$.
\end{enumerate}
\end{fact}

\pf
Setting $\ue = (d_1,1,d_1)$ and $i^\ast = 3$, 
Lemma~\ref{translation.lemma} implies that 
$d_1 v(x_2) = v(x_1^{d_1}) + d_1 v(x_2) \geq n$.
Therefore $v(x_2) \geq n/d_1 > 1$, so the first statement 
holds for $i = 2$.  
Using $\ue = (d_1,d_1,1)$ and $i^\ast = 2$ establishes 
the statement for $i = 3$.
To establish the second statement, use
$\ue = (1,d_2,d_2), i^\ast = 3$.
The remaining two statements follow from
{\bf\ref{zeta.sub.s}}.
\eop

\smallskip
Set $\displaystyle \zeta^t(d) = 
       \frac1d\left(1 + \sum_{m|d,m<d} \phi(d/m) p^{-\max(1,n-mt)}\right)$.

Note that 
$$\zeta^t(d) = \zeta_{n-t,n-2t, \ldots}(d).$$

\begin{lemma}
\label{LemmaX}
If $j \neq k$ and $\sum_{i \neq j,k} v(x_i) \leq t$, then $\kappa(x_j) \leq \zeta^t(d_j)$
and $d_j \geq n/t$.
\end{lemma}

\pf
Without loss, $j = 1$ and $k = 2$.
From Lemma~\ref{translation.lemma} with $\ue = (m,m,1,\ldots,1)$
and $i^\ast = 2$ 
we have $v(x_1^m) \geq n - mt$.
The total contribution of the $\phi(d_1/m)$ generators of $\lan x_1^m \ran$ to $\kappa(x_1)$
is therefore at most $\phi({d_1}/m) \cdot \frac1{d_1} \cdot p^{-\max(1,n-mt)}$.
This implies the inequality for $\kappa(x_j)$.  Since $v(x_1^{d_1}) =  0$, it
also follows that $d_1 \geq n/t$.
\eop

\begin{lemma}
\label{LemmaX4}
If $j, k, l$  are distinct, $d_k = d_l = 2$, and 
$\sum_{i \neq j,k,l} v(x_i) \leq t$, then $\kappa(x_j) \leq \zeta^{2t}(d_j)$
and $d_j \geq n/2t$.
\end{lemma}

\pf
Argue as in the proof of Lemma~\ref{LemmaX}.
Assume $j = 1$, $k = 2$, $l = 3$, and use Lemma~\ref{translation.lemma}
with $\ue = (m,2,2,1,\dots,1)$ and $i^\ast = 1$ to get
$v(x_1^m) \geq n - 2mt$.  
\eop

\begin{lemma}
\label{LemmaY}
Suppose $\ud = (2,d_2,d_3)$ and $v(x_2^2) = v$. 
\begin{enumerate}
\item
$\kappa(x_3) \leq \zeta^v(d_3)$ and $d_3 \geq n/v + 1$.
\item
If $p = 2$ then 
$\kappa(x_3) \leq \zeta^{v/2}(d_3)$ and $d_3 \geq 2n/v$.
\end{enumerate}
\end{lemma}

\pf Using $\ue = (2,2,k)$, $i^\ast = 3$, in
Lemma~\ref{translation.lemma},  
we have $v(x_3^k) \geq n - kv$ in general, and
$v(x_3^k) \geq n - kv/2$ when $p = 2$.  
Using $\ue = (2,2,d_3)$, $i^\ast = 2$,
we have $(d_3-1)v \geq n$ in general 
and $d_3 v \geq 2n$ when $p=2$.
\eop

\begin{lemma}
\label{Scott3}
If $r = 3$ and $i \neq j$, then $d_i v(x_j) \geq n$.
In particular, $\kappa(x_j) \leq \zeta_{\ceil{n/d_i}}(d_j)$, 
where $\ceil{x}$ is the smallest integer not less than $x$.
\end{lemma}

\pf Without loss, $i = 1$ and $j = 2$.
The first statement follows from 
Lemma~\ref{translation.lemma} with $\ue - (d_1,1,d_1)$ and $i^\ast = 3$.
The second statement follows from the first.
\eop

\begin{lemma}
\label{23d}
Assume that $\ud = (2,3,d)$.
If $p$ is odd, set 
$s_2 = \ceil{n/2}$,
$s_3 = \ceil{n/3}$,
$s_4 = \ceil{n/5}$, and
$s_5 = \ceil{n/11}$.
If $p = 2$, set 
$s_2 = \ceil{2n/3}$, 
$s_3 = \ceil{n/2}$, 
$s_4 = \ceil{n/3}$, and
$s_5 = \ceil{n/6}$.
Then
$$v(x_3^k) \geq s_k, \ \ d = 2,3,4,5.$$
In particular $\kappa(x_3) \leq \zeta^\ast_{s_2,s_2,s_3,s_4,s_5}(d)$.
\end{lemma}

\pf
Lemma~\ref{translation.lemma} with $\ue = (2,3,e)$ and $e = 2,3,4,5$,
with $i^\ast = 3$ shows that $(C_3(\ue) -1)v(x_3^e) \geq n$ in general 
and $C_3(\ue) v(x_3^e) \geq 2n$ when $p = 2$.
We have $C_3(\ue) = 3,4,6,12$ in the respective situations, and the 
result follows immediately. 
\eop

\begin{lemma}
\label{24d}
Assume that $\ud = (2,4,d)$.
If $p$ is odd, set 
$s_2 = \ceil{n/3}$ and
$s_3 = \ceil{n/7}$.
If $p = 2$, set 
$s_2 = \ceil{n/2}$ and
$s_3 = \ceil{n/4}$.
Then
$$v(x_3^k) \geq s_k, \ \ d = 2,3.$$
In particular $\kappa(x_3) \leq \zeta^\ast_{s_2,s_2,s_3}(d)$.
\end{lemma}

\pf
Use Lemma~\ref{translation.lemma} with $\ue = (2,4,e)$ and $e = 2,3$,
with $i^\ast = 3$ for the general case.
We have $C_3(\ue) = 4,8$ in the respective situations.
\eop

\begin{lemma}
\label{zeta5*}
Suppose $p=2$, $n \geq 14$, $r=3$, and $\{i,j,k\} =  \{1,2,3\}$.
Then 
\begin{enumerate}
\item 
\label{part1}
$v(x_i^2) + v(x_j^2) \geq 28/d_k$.
\item 
\label{part2}
If $d_i = d_j = 3$, then $v(x_k^2) \geq 5$.
\item 
\label{part3}
If $d_i = 3$ and $d_j = 4$, then $v(x_k^2) \geq 3$.
\end{enumerate}
\end{lemma}

\pf
Without loss, $i = 1$, $j = 2$, and $k = 3$.
Use Lemma~\ref{translation.lemma} with $\ue = (2,2,d_3)$, $(3,3,2)$, and
$(3,4,2)$, respectively.
\eop

\subsection{Initial reductions}

The proof of Theorem~\ref{basic.result} uses 
routine, but extensive, calculations based on the 
results of the previous subsections.
We have verified these calculations using GAP4 \cite{GAP4}.

Assume that 
\begin{enumerate}
\item $G$ is a classical group with natural module $V$ and $\F_p$ dimension $n$.
\item $V$ contains at least $10^4$ points.
\item  $\ux$ is a normalized generating $r$-tuple 
       for $G$ in a primitive action.
\item Every power of every element of $\ux$ satisfies 
  Grassmann Condition $1/100$.
\item $g(\ux) \leq 2$.
\end{enumerate}
To prove Theorem~\ref{basic.result} 
it suffices to  show that the characteristic of $V$, the dimension 
of $V$ over its prime field, and the signature of $\ux$ are given 
in Table~\ref{pnd.table}.

Unless stated otherwise, we assume that $d_1 \leq d_2 \leq \dots \leq d_r$
and that $v(x_i) \leq v(x_{i+1})$ whenever $d_i = d_{i+1}$.  
Also recall that $\epsilon_0$ and $A(d)$ were defined just before Fact {\bf \ref{kappa.inequality}}.

We have $\eps_0 < 2 \cdot 10^{-4}$ and $\eps < 10^{-2}$.
Combining Fact {\bf \ref{kappa.inequality}} 
with the inequality $\kappa(x_i) \leq \frac1{d_i} + \frac1p - \frac1{d_ip}$,
we have the following inequalities.

\begin{fact}
\label{Ap.fact}
$A(\ud) > (.99A(\ud) - 2.0002)p$.
Consequently
\begin{enumerate}
\item
\label{p.inequality}
$\displaystyle p <
\frac{A(\ud)}{.99A(\ud) - 2.0002}$
\item
\label{A.inequality}
$\displaystyle 
A(\ud) < \frac{2.0002p}{.99p - 1}$
\end{enumerate}
\end{fact}

\begin{lemma} 
\label{n.bounds}
$n \geq 3$.
\begin{enumerate}
\item If $p \leq 97$ then $n \geq 4$.
\item If $p \leq 19$ then $n \geq 5$.
\item If $p = 7$ then $n \geq 6$.
\item If $p = 5$ then $n \geq 7$.
\item If $p = 3$ then $n \geq 10$.
\item If $p = 2$ then $n \geq 14$.
\end{enumerate}
\end{lemma}

\pf The enumerated statements are immediate consequences
of the inequality $(p^n-1)/(p-1) \geq 10000$.

If $n = 2$, then $F^\ast(G) \cong L_2(p)$, and $F(x) \leq 2$ for all
$x \in G^\sharp$.  It follows that $f(x) \leq 1/5000$ for all $x \in G^\sharp$,
so equation (\ref{CF.consequence}) cannot hold for $g \leq 2$.  
\eop

Set $S = S(\ud) = r-2 - .01A(\ud) - .0002$, the right hand side of the inequality
in statement {\bf \ref{kappa.inequality}}.
For $i = 1, \ldots, r$, set $\kappa_i = \kappa(x_i)$.
Set $\Sigma = \sum \kappa_i$.  
Then $\Sigma > S$ by {\bf\ref{kappa.inequality}} and assumptions 
on $\ux$.

\begin{lemma}
\label{p.prop}
\begin{enumerate}
\item
If $p \geq 17$, then $r = 3$.
\item
If $p \geq 7$, then $r \leq 4$.
\item
If $p = 7$, then $r \leq 4$ and $S \geq (r-3) + .9761$.
\item
If $p = 5$, then $r \leq 5$ and $S \geq (r-3) + .9744$.
\item
If $p = 3$, then $r \leq 6$ and $S \geq (r-3) + .9693$.
\item
If $p = 2$, then $r \leq 8$ and $S \geq (r-3) + .9589$.
\end{enumerate}
\end{lemma}

\pf
Since 
$\zeta$ is a decreasing function, we have
$\zeta(d) \leq \zeta(2) = (p+1)/2p$, so $\kappa(x_i) \leq (p+1)/2p$
for all $i$.  Therefore $r \cdot \frac{p+1}{2p} > r - 2 - .01A(\ud) - .0002 > .99r -2.0002$,
whence
$$r < \frac{4.0004p}{.98p - 1}.$$
All assertions about $r$, except the first, follow from this.

If $r = 4$, then $A(\ud) \geq 13/6$, so $p < 17$ by {\bf\ref{Ap.fact}.\ref{p.inequality}}.

The statements concerning $S$ follow from {\bf \ref{Ap.fact}.\ref{A.inequality} }.
\eop

\begin{lemma}
\label{S.bounds}
\begin{enumerate}
\item If $r = 3$, then $S \geq .9698$.
\item If $r=3$ and $d_1 = 2$ then $S \geq .9748$.
\item If $r=3$, $d_1 = 2$, and $d_2 = 3$ then $S \geq .9781$.
\item If $\ud = (2,3,7)$, then $S \geq .9795$.
\end{enumerate}
\end{lemma}

\pf
These statements follow from straightforward computations.
\eop

\subsection{Completion of the Proof}


\begin{lemma}
$n \geq 4$.
\end{lemma}

\pf
Suppose $n = 3$.  Then $\Omega$ is the set of points in the natural
module for $F^\ast(G) \cong L_3(p)$.  We have $N = p^2 + p + 1$.
By Lemma~\ref{n.bounds}, $p > 100$, so $A(\ud) < 2.02$ by
{\bf\ref{Ap.fact}.\ref{A.inequality}}.  It follows that $\ud = (2,3,7)$.

Since $x_1$ is an involution in $G$, we have $\Fix(x_1) \leq p+2$,
and $\Ind(x_1) \geq \frac12(p^2-1)$.
By Lemma~\ref{Scott3},
$v(x_i) \geq 2$, for $i = 2,3$.
This implies that $\Fix(x_i) \leq 3$, $i = 2, 3$, whence
$\Ind(x_i) \geq (d_i-1)/d_i \cdot (p^2 + p -2)$.
It follows from the Riemann-Hurwitz equation that $g > 2$, a contradiction.
\eop

\begin{lemma}
\label{prop.lt.23}
$p \leq 19$
\end{lemma}

\pf
Suppose $p \geq 23$.
Then $A(\ud) \leq 2.0002p/(.99p-1) < 2.114$
by {\bf \ref{Ap.fact}.\ref{A.inequality} }.
This implies that $\ud$ is one of the following:
$(2,3,d)$, $(2,4,\leq7)$, $(2,5,5)$, or $(3,3,4)$.
Also, $S > .9787$ by {\bf\ref{kappa.inequality}}.

If $\ud = (2,3,d)$, $d \geq 8$,
then {\bf \ref{zetas}} implies that 
$\sum \kappa(x_i) \leq \zeta_2(2) + \zeta_2(3) + \zeta_2(d)$.
Since $\phi(d) \geq 4$, it follows that
$\zeta_2(d) \leq \frac1d(1 + (d-5)/p + 4/p^2) =
\frac1p + (1 + \frac4{p^2} + \frac5p)\cdot\frac1d 
\leq \zeta_2(8) < .1423$, 
whence $\sum\zeta_2(d_i) \leq .9778$, a contradiction.

In the remaining six cases, we have $\kappa_i \leq \zeta_2(d_i), i = 2,3$ 
and $\kappa_1 \leq \zeta(d_1)$ in all cases, and
$\kappa_1 \leq \zeta_2(2)$ in the $(2,3,7)$ case.
By inspection, either $\Sigma < S$ or $p=23$ and $\ud = (2,4,5)$ or
$(2,3,7)$.

Suppose $\ud = (2,4,5)$.  Then $v(x_3) \geq \mustar(5,23) = 4$.
If $v(x_1) \geq 2$, then $\sum \kappa_i < .9628$, so we must have
$v(x_1) =1$.  Therefore $v(x_2) \geq n - v(x_1) \geq 3$.  
Furthermore, $n = 4$, by
Lemma~\ref{Scott3}.
This implies that
$v(x_2^2) \geq 2$ since every involution $t$  in $PGL(4,23)$
with $v(t) = 1$ is not a square in that group.
It follows that $\sum \kappa_i < .974$, a contradiction.

We must have $\ud = (2,3,7)$, whence $v(x_3) \geq \mustar(7,23) =3$,
and $\kappa_3 \leq \zeta_3(7)$.  This implies that $\sum \kappa_i < .9786$,
which is not so.
\eop

\begin{prop}
\label{large.p}
If $p > 7$ then $p = 11$, $\ud = (2,3,7)$, and $n = 5$ or $6$.
\end{prop}

\pf
By Lemmas~\ref{n.bounds} and~\ref{prop.lt.23}, $n \geq 5$.
Suppose $p > 7$.
Then $p \geq 11$, and for purposes of estimation with $\zeta(d)$ 
and $\zeta_k(d)$ we
may assume that $p = 11$.

Since $A(\ud) \leq 2.2246$ by {\bf\ref{Ap.fact}}, we have $S \geq (r - 3)
+ .9775$.

If $r > 3$, then $\ud = (2,2,2,3)$ by the condition on $A(\ud)$.
Since $\sum_{i \neq j} v(x_i) \geq n \geq 5$ for $j = 3, 4$, we
have $v(x_3) > 1$ and either $v(x_2) > 1$ or $v(x_4) > 1$.
Therefore $\sum \kappa_i \leq \max(2 \zeta(2) + \zeta_2(2) + \zeta_2(3),
\zeta(2) + 2\zeta_2(2) + \zeta(3))< 1.95$, a contradiction.

Thus $r = 3$.  Since $\zeta(d_1) \geq S/3 > \zeta(4)$, it
follows that $d_1 = 2$ or $3$.

Suppose $d_1 = 3$.  Then $\zeta(d_2) > (S-\zeta(3))/2 > \zeta(5)$,
so $d_2 \leq 4$, and $\kappa_1 \leq \zeta_2(3)$ by {\bf \ref{zetas}}.
Since $\zeta_2(4) < \zeta_2(3)$, this implies that $\kappa_3 >
S - 2 \zeta_2(3) > \zeta_2(4) > \zeta(5)$, whence $d_3 = 3$,
which is impossible because $\ud \neq (3,3,3)$.  This shows that 
$d_1 = 2$.

Since $\kappa_3 \leq \zeta(d_3) \leq \zeta(d_2)$ and 
$\kappa_2 \leq \zeta(d_2)$, we must have
$\zeta(d_2) > (S- \zeta(2))/2 > \zeta(8)$ so $d_2 \leq 7$.
If $d_2 = 5, 6$, or $7$, then $\kappa_2 \leq 
\max_{5 \leq d \leq 7}(\zeta_2(d)) \leq \zeta_2(6)$.
[Recall that $p=11$ for the purpose of calculation.]
Since $\zeta(8) < \zeta_2(6)$ and $\zeta(d) < \zeta(8)$ for $d > 8$, 
we have $\kappa_3 \leq \zeta_2(6)$
and $\sum \kappa_i \leq \zeta(2) + 2 \zeta_2(6) < S$.
Therefore $d_2 \leq 4$.

Suppose $d_2 = 4$.
Then $\kappa_3 \geq S - \zeta_2(2) - \zeta_3(4) > .2002 > \zeta_3(d)$
for $d > 6$, so $d_3 \leq 6$.  If $d_3 = 5$, then $A(\ud) = 2.05$,
so $S \geq .9793$ and $\Sigma \leq \zeta_2(2) + \zeta_3(4) + \zeta_3(5) < .9781$.
It follows that $d_3 = 6$.
From Lemma~\ref{translation.lemma} with $\ue = (2,2,2)$,
either $v(x_2^2) > 1$ or $v(x_3^2) > 1$.
If $v(x_2^2) > 1$, then $\kappa_2 \leq \zeta_{3,2}(4)$.
If $v(x_3^2) > 1$, then $\kappa_3 \leq \zeta_{3,2}(6)$.
In either case, $\Sigma < .97 < S$.  

Suppose $d_2 = 3$.
Then $S \geq .9781$ and $\kappa_1 + \kappa_2 \leq \zeta_2(2) + \zeta_3(3)
< .8426$, so $\kappa_3 > .1355$.
If $d \geq 21$, then $\zeta(d) < \zeta(21) < .135$.
Therefore $d_3 \leq 20$.  By inspection, if
$d_3 = 9$ or $d_3 \geq 11$, then $\zeta_3(d_3) < .137$, and the inequality
cannot hold.  
Therefore $d_3$ is one of $7, 8$, or $10$.
If $d_3 = 8$, or $10$, then 
$\kappa_3 \leq \zeta_{3,3}(d_3)$
by Lemma~\ref{23d} and $\sum \kappa_i < S$. 

Therefore $d_3 = 7$, so $S  \geq .9795$ and
the condition $\zeta_2(2) + \zeta_3(3) + \zeta_3(7) \geq S$
implies that $p = 11$ or $13$.
If $p = 13$, then $v(x_3)$ is necessarily even, so $\kappa_3 \leq \zeta_4(7)$
and $\sum \kappa_i < S$.
Therefore $p = 11$.
It follows that $\kappa_2$ is even, so $\kappa_2 \leq \zeta_4(3)$.
If $v(x_1) > 2$, then $\sum \kappa_i \leq \zeta_3(2) + \zeta_4(3) + \zeta_3(7) < S$.
Therefore $v(x_1) = 2$ and $n = 5$ or $6$.
\eop

\begin{prop}
\label{p=7.prop}
If $p = 7$, then $n = 6$ and $\ud = (2,3,7)$.
\end{prop}

\pf
By Lemma~\ref{p.prop}
$n \geq 6$, $r \leq 4$, and 
$S \geq (r-3) + .9761$.

Suppose $r = 4$.
If $v(x_1) + v(x_2) = 2$, then $d_j \geq 3$, $j > 2$, and $\sum \kappa_i \leq 2 \zeta(2)
+ 2 \zeta^2(3) < 1.9$ by Lemma~\ref{LemmaX}.
Therefore $v(x_1) + v(x_2) \geq 3$, and in fact $v(x_i) \geq 2$ for at least $3$
choices of $i$.  It follows from inspection
of values of $\zeta(d)$ and $\zeta_2(d)$ that $\sum \kappa_i < S$, a contradiction.

Therefore $r = 3$.  If $v(x_1) = 1$, then $\kappa_2, \kappa_3 \leq \zeta^1(d) < .168$
by Lemma~\ref{LemmaX}.   
Since $\kappa_1 \leq \zeta(2) < .572$, we have $\sum \kappa_i
< S$, a contradiction.  Therefore $v(x_i) \geq 2$ and $\kappa_i \leq \zeta_2(d_i)$
for all $i$.  Since $\zeta_2(d) < .3$ for $d > 3$, we have $d_1 \leq 3$.

Suppose $d_1 = 3$.  Then, by inspection of  $\zeta_2(d)$, $d \geq 3$, we have
$\ud = (3,3,4)$.  Either $v(x_1) = 2$, in which case $\sum \kappa_i \leq
\zeta_2(3) + \zeta_4(3) + \zeta_4(4)$, or $v(x_1) \geq 3$, in which case
$\sum \kappa_i \leq 2 \zeta_3(3) + \zeta_2(4)$.  In either case, $\sum \kappa_i
< .97$, a contradiction.  We conclude that $d_1 = 2$.

We have $\kappa_2 + \kappa_3 \geq S - \zeta_2(2) \geq .465$.
Also $\kappa_i \leq \zeta_3(d_i), i > 1$ by Lemma~\ref{Scott3}.
By inspection, $\zeta_3(d) < .2$ for $d > 6$, so $d_2 \leq 6$.

Suppose $v(x_1) = 2$.  Then $\kappa_j \leq \zeta^2(d_j), j \geq 2$, whence
$d_2 \leq 4$ because $\zeta^2(d) < .21$ for $d > 4$.
If $d_2 = 4$, then $d_3 \geq 5$ because $A(\ud) > 2$, so
$\sum \kappa_i \leq \zeta_2(2) + \zeta^2(4) + \zeta^2(5) < .97$, a contradiction.
Therefore $d_2 = 3$ and $d_3 \geq 7$, so 
$S \geq .9781$, and $\kappa_2 + \kappa_3 \geq .4678$,
whence $\kappa_3 \geq .4678 - \zeta^2(3) > .1341$.  
By inspection,
$d_3 \in \{7,8,9,12\}$.
By Lemma~\ref{23d}, $\kappa_3 < \zeta_{4,3,2,2}(d_3)$, and
we conclude that $\ud = (2,3,7)$.
Note that $n = 6$ by Lemma~\ref{Scott3}.

We may assume henceforth that $v(x_1) \geq 3$, so $\kappa_1 \leq .5015$ and 
$\kappa_2 + \kappa_3 > .4747$. 

If $d_2 = 6$, then $d_3 = 6$ by inspection
of the values of $\zeta_3(d)$, $d \geq 6$.  
From Lemma~\ref{translation.lemma} with $\ue = (2,2,2)$
we have $v(x_j^2) > 1$ for some $j > 1$, so 
$\kappa_2 + \kappa_3 \leq \zeta_3(6) + \zeta_{3,2}(6) < S -\kappa_1$.
This implies that $d_2 < 6$.  

By inspection, $d_2 \neq 5$.  
If $d_2 = 4$, then $\kappa_2 \leq \zeta_3(4) < .2872$, so $\kappa_3 > .1875$.  
This implies
that $d_3 \leq 6$.  From Lemma~\ref{translation.lemma} with $\ue = (2,2,d_3)$
we have
$v(x_2^2) > 1$, so 
$\kappa_2 \leq \zeta_{3,2}(4)
< .257$.  
When $d_3 = 6$, the same argument shows that $\kappa_3 \leq 
\zeta_{3,2}(6)< .2$.  
In each case, $\sum \kappa_i < S$.

So $d_2 \neq 4$, and we have $d_2 = 3$.
Also, $\kappa_1 + \kappa_2 \leq \zeta_3(2) + \zeta_3(3) < .8368$.
so $\kappa_3 \geq S - \kappa_1 - \kappa_2 > .1413$.
By inspection of $\zeta_3(d)$, we have $d_3 \leq 18$.
By Lemma~\ref{23d}, $\kappa_3 \leq \zeta_{3,3,2,2}(d_3)$, so by inspection $d_3 = 7$.
If $n > 6$, then $v(x_1) \geq 3$ and $v(x_j) \geq 4$, $j > 1$,
so $\sum \kappa_i \leq \zeta_3(2) + \zeta_4(3) + \zeta_4(7) < .9781 < S$.
Therefore $n = 6$.
\eop

\begin{prop}
\label{p=5.prop}
If $p = 5$, then $\ud = (2,3,7)$, $n = 7$, $8$, or $9$, $v(x_1) = 3$, and $v(x_3) = 6$.
\end{prop}

\pf
By Lemma~\ref{p.prop},
$n \geq 7$, $r \leq 5$, and 
$S \geq (r-3) + .9744$.

If $r = 5$, then $\sum \kappa_i  \leq 3 \zeta(2) + 2 \zeta_2(2) < S$ because
$v(x_i) > 1$ for at least two choices of $i$.  Therefore $r \leq 4$.

Suppose $r = 4$.  
If $v(x_1) + v(x_2) \leq 3$, then 
Lemma~\ref{LemmaX} implies that  $d_i \geq 7/3 > 2$ for $i = 3, 4$, and $\kappa_i
\leq \zeta^3(d_i) \leq \zeta^3(3) = .3344$. 
Since $\kappa_1 + \kappa_2 \leq 2 \zeta(2) = 1.2$, it follows that
$\Sigma < S$.
Therefore $v(x_1) + v(x_2) \geq 4$.  
Moreover, $v(x_i) + v(x_j) \geq 4$ whenever $i \neq j$.
If $v(x_1) = 1$, then $\sum \kappa_i \leq \zeta(2) + 2 \zeta_3(2) + \zeta_3(3) < 1.95$.
If $v(x_1) = 2$, then $\sum \kappa_i \leq 3\zeta_2(2) + \zeta_2(3) = 1.92$.
Therefore $v(x_1) \geq 3$.
If $d_3 > 2$, then 
we have$\sum \kappa_i < 2 \zeta_3(2) + 2 \zeta(3) < 1.95$, 
noting that $\kappa_1, \kappa_2 \leq \zeta_3(2)$ since
$\zeta(3) < \zeta_3(2)$.  
So $d_3 = 2$ and $\kappa_1 + \kappa_2
+ \kappa_3 \leq 3 \zeta_3(2) = 1.512$.  
From Lemma~\ref{translation.lemma} with $\ue = (2,2,2,1)$
we have
$v(x_4) \geq 3$, so $\kappa_4 \leq \zeta_3(d)$ and $\sum \kappa_i \leq 3 \zeta_3(2)
+ \zeta_3(d) < 1.9$.  We conclude that $r \neq 4$.  Thus, $r = 3$.

\medskip
If $v(x_1) = 1$, then Lemma~\ref{LemmaX} shows that $d_2 \geq 7$ and
$\kappa_i \leq \zeta^1(d_i)$, $i = 2, 3$.  So $\sum \kappa_i \leq \zeta(2) + 2 \zeta^1(d)
< .9$.  Therefore $v(x_1) \geq 2$, and in fact $\kappa_i \leq \zeta_2(d_i)$ for
all $i$.  
Since $\zetao(d) \leq .29$ for $d \geq 12$ by Lemma~\ref{Lemma5} and 
$\zeta_2(d) \leq .32$ for $4 \leq d \leq 11$ by inspection
it follows that $d_1 \leq 3$, whence $\kappa_i \leq \zeta_3(d_i)$,
$i = 2,3$ by Lemma~\ref{Scott3}.

Suppose $d_1 = 3$.  Then $\kappa_1 \leq \zeta_2(3) = .36$.
If $d_2 \geq 4$, then $\kappa_i \leq \zeta_3(d) \leq.304$ for $i > 1$, and
$\sum \kappa_i \leq .968$.
Therefore $d_2 = 3$, so $v(x_1) \geq 3$ and $\kappa_1 + \kappa_2 \leq 2 \zeta_3(3)
< .6774$.  If $d_3 > 4$, then $\kappa_3 \leq \zeta_3(d_3) \leq .27$, so $d_3 = 4$.
From Lemma~\ref{translation.lemma} with $\ue = (3,3,2)$
we have $v(x_3^2) \geq 2$ so
$\kappa_3 \leq \zeta_{3,2}(4) = .264 
< S - \kappa_1 - \kappa_2$.
We conclude that $d_1 = 2$.

We have shown that $v(x_1) > 1$.
If $d_3 > 23$, then $\kappa_3 < \zetao(d_3) < .165$.
Suppose $v(x_1) = 2$.  Then $d_2 \geq 4$ and $\kappa_i \leq \zeta^2(d_i)$, $i= 2,3$.
Since $\kappa_1 \leq .52$ and $\zeta^2(d) \leq .203$ for $d > 4$,
we must have $d_2 = 4$, whence $d_3 > 4$.  If $d_3 \neq 6$, then
$\sum \kappa_i \leq \zeta_2(2) + \zeta^2(4) + \zeta^2(5) < .973 < S$.
Therefore $d_3 = 6$ and  $A(\ud) <  2.09$, so $S > .9789$.
We have $\sum \kappa_i \leq \zeta_2(2) + \zeta^2(4) + \zeta^2(6) \leq .975$,
a contradiction.  This shows that $v(x_1) > 2$.

We have $\kappa_1 \leq \zeta_3(2) = .504$ 
so $\kappa_2 + \kappa_3 > .47$.
Also, $\kappa_i \leq \zeta_4(d_i), i = 1,2$, 
Since 
$\zetao(d) < .2$ for $d \geq 20$ and
$\zeta_4(d) < .21$ for $6 < d < 20$, we have $d_2 \leq 6$. 
From Lemma~\ref{translation.lemma} with $\ue = (2,d_2,2)$ and $i^\ast = 3$
it follows that
$v(x_3^2) \geq 7/(d_2-1) > 1$.
This implies that $\kappa_3 \leq \zeta_{4,2}(d_3)$.
If $d_2 = 6$, then $\kappa_2 \leq .268$.
We have $d_3 = 6$ as otherwise $\kappa_3  \leq
\min(\zetao(d_3),\zeta_{4,2}(d_3)) < .174$,
so $\kappa_2, \kappa_3 \leq \zeta_{4,2}(6) < .214$, 
and $\Sigma < S$. 
Therefore $d_2 < 6$.
If $d_2 = 5$, then $\kappa_2  + \kappa_3 \leq \zeta_4(5) + \zeta_{4,2}(6) \leq .47$.
If $d_2 = 4$, then $\kappa_2 \leq \zeta_4(4) = .3008$.
By Lemma~\ref{24d},
$\kappa_3
\leq \zeta_{4,3}(d_3) $. 
If $d_3 > 23$, then $\kappa_3 < \zetao(d_3) < .165$.
It follows from inspection that 
$\zeta_{4,3}(d)\leq .214$ for $6 < d < 24$.
Therefore $d_3 \leq 6$, 
whence $v(x_2^2) \geq 2$ and $\kappa_2 \leq
\zeta_{4,2}(4) = .2608$.
If $d_3 = 5$, then $\sum\kappa_i < S$.
If $d_3 = 6$, then $\kappa_3 \leq .2032$, and $\sum \kappa_i < S$.
It follows from this paragraph that $d_2 \neq 4$.  Therefore $d_2 = 3$.

We have $\kappa_1 + \kappa_2 \leq \zeta_3(2) + \zeta_4(3) = .8384$.
By Lemma~\ref{p.prop},
$S \geq .9781$, so $\kappa_3 \geq .1397$.
Since $\kappa_3 < 3/d_3 +.04$ we may assume that $d_3 \leq 30$.
By Lemma~\ref{23d}, $\kappa_3 \leq \zeta_{4,4,3,2}(d_3)$.
By inspection, $d_3 = 7$.

If $v(x_1) \geq 4$, then $\sum \kappa_i \leq 
\zeta_4(2) + \zeta_4(3) + \zeta_4(7) < S$.
So $v(x_1) = 3$ and $n \leq d_2 v(x_1) = 9$.
\eop

\begin{prop}
\label{p=3.prop}
If $p = 3$, then either
\begin{enumerate}
\item
$\ud = (2,3,7)$, $n = 12$, $v(x_1) =  4$, $v(x_2) = 8$, and $v(x_3) = 12$ or 
\item
$\ud = (2,3,8)$, $n = 10$, $v(x_1) = 4$, $v(x_2) = 6$, and $v(x_3^4) = 2$.
\end{enumerate}
\end{prop}

\pf
By Lemma~\ref{p.prop},
$n \geq 10$, $r \leq 6$, and 
$S \geq (r-3) + .9693$.

We note that $\zetao(d) < .11$ for $d > 42$ by Lemma~\ref{Lemma5}
and $\zetao(d) < .11$ by direct computation for $24 < d \leq 42$.
Also, $\zetao(d) < .2$ for $d > 12$.
Thus, statements bounding $\kappa_i$ with weaker bounds need 
only be verified for a finite number of possible values of $d_i$.
We shall use this implicitly in the following argument.

Since $n > r$, we have $\kappa_i \leq \zeta_2(d_i)$ for at least two choices of $i$.
If $r= 6$, then $\sum \kappa_i \leq 4 \zeta(2) + 2 \zeta_2(2) < 3.8$, a contradiction,
so $r \leq 5$.

Suppose $r= 5$.  
If $v(x_1)+ v(x_2) + v(x_3) = 3$, then $d_i \geq 4$ and
$\kappa_i \leq \zeta^3(d_i) < .3$ for $i = 4, 5$
by Lemma~\ref{LemmaX}.  
If $v(x_1)+ v(x_2) + v(x_3) = 4$, then $d_i \geq 3$
and $\kappa_i \leq \zeta^4(d_i) < .35$ for $i = 4,5$
Since $\kappa_1 + \kappa_2 + \kappa_3 \leq 3 \zeta(2) = 2$ we have
$\Sigma < S$ in this case.
Therefore $v(x_i) + v(x_j) + v(x_k) \geq 5$ for any choice of distinct $i,j,k$.
If $v(x_i) = 1$ for two values of $i$, then $v(x_i) \geq 3$ for
three values and $\sum \kappa_i \leq 2\zeta(2) + 3 \zeta_3(2) < 2.9$.
Therefore $v(x_i) = 1$ for at most one value of $i$, and $\sum \kappa_i 
\leq \zeta(2) + 4\zeta_2(2) < 2.9$.  We conclude that $r \leq 4$.

Suppose $r= 4$.
If $v(x_1) + v(x_2) = 2,3,4$, respectively, then 
$\kappa_1 + \kappa_2$ is respectively  at most $1.3334$, $1.2223$, $1.1852$,
while Lemma~\ref{LemmaX} implies that for $i = 3$ or $4$,
$\kappa_i \geq \zeta^2(d_i)$ and $d_i \geq 5$,
$\kappa_i \geq \zeta^3(d_i)$ and $d_i \geq 4$,
$\kappa_i \geq \zeta^4(d_i)$ and $d_i \geq 3$,
in the respective cases.
By inspection, $\kappa_3 + \kappa_4$ is respectively
at most $.401$, $.511$, $.67$, whence $\sum \kappa_i < S$.
It follows that $v(x_1) + v(x_2) \geq 5$.
Since the same is true of $v(x_i) +  v(x_j)$, $i \neq j$, it follows that
$v(x_i) \geq 3$ for at least $3$ choices of $i$.
Since $\zeta(2) < .67$, $\zeta_3(2) < .52$, and $\zeta(d) < .56$, 
$\zeta_3(d) < .36$ when $d > 2$, we have $d_3 = 2$, else $\sum \kappa_i < 1.96 < S$.
Set $v = v(x_1)$.  If $v= 1$, then $\kappa_1 + \kappa_2 + \kappa_3 
\leq \zeta(2) + 2 \zeta_4(2) < 1.68$ and, by Lemma~\ref{LemmaX4},
$\kappa_4 \leq \zeta^2(d_4)$ where $d_4 \geq 5$, so $\kappa_4 < .21$.
If $v= 2$, then $\kappa_1 + \kappa_2 + \kappa_3 
\leq \zeta_2(2) + 2 \zeta_3(2) < 1.6$ and, by Lemma~\ref{LemmaY}, and
inspection of $\zeta^2$ values,
$\kappa_4 \leq \zeta^4(d_4)$ where $d_4 \geq 4$, so $\kappa_4 < .28$.
If $v > 2$, then $\kappa_1 + \kappa_2 + \kappa_3 \leq 3 \zeta_3(2) < 1.56$.
From Lemma~\ref{translation.lemma} with $\ue = (2,2,2,1)$ and $i^\ast = 4$
we have $v(x_4) \geq 4$
and $\kappa_4 \leq \zeta_4(3) < .35$.
In all cases, $\sum \kappa_i < S$.
Therefore $r \neq 4$.

We have $r=3$.  By inspection, $d > 6$ implies $\zetao(d) \leq .25$.
Therefore $d_1 \leq 6$ and $A(\ud) \leq 3 \cdot 5/6 < 2.84$, so $S > .9714$.
By inspection, $\kappa_1 \leq \zeta(2) < .67$.

If $v(x_1) = 1$, then, by Lemma~\ref{LemmaX}, $d_i \geq 10$ and $\kappa_i \leq
\zeta^1(d_i) < .11$, $i = 2,3$.
If $v(x_1) = 2$, then $\kappa_1 \leq \zeta_2(2) < .556$. 
Also, by Lemma~\ref{LemmaX}, $d_i \geq 5$ 
and $\kappa_i \leq \zeta^2(d_i) < .201$, $i = 2,3$.
It follows that $v(x_1) > 2$, so $\kappa_i \leq \zeta_3(d_i)$ for all $i$.

Suppose $d_1 \geq 4$.  Then $\kappa_i \leq \zeta_3(d_i) \leq \zeta_3(4) < .3519$
for all $i$, so $\sum_{j \neq i} \kappa_j \geq S - \zeta_3(4) > .619$,
$i = 1,2,3$.  
If $v(x_i) = 3$ for some $i$,
then Lemma~\ref{LemmaX} shows that $\kappa_j \leq \zeta^3(d_j) \leq .254$
for $j \neq i$.  
If $v(x_i) = 4$ for some $i$, then $\kappa_i \leq 
\zeta_4(d_i) \leq \zeta_4(4) < .34$ and $\kappa_j \leq \zeta^4(d_j) < .28$,
$j \neq i$.  
It follows that $v(x_i) \geq 5$ for all $i$, so $\kappa_i \leq \zeta_5(4) < .336$.
Since $\zetao(d_3) \leq .25$ for all $d > 4$ with $d \neq 6$ we conclude that
$d_i = 4$ or $6$ for all $i$.  
From  Lemma~\ref{translation.lemma} with $\ue = (2,2,2)$
we have
$v(x_i^2) \geq 2$ for some $i$.  
Therefore $\kappa_i \leq \zeta_{5,2}(d_i) < .28$ for some $i$.
Since $\kappa_j \leq \zeta_5(d_j) < .34$ for all $j$ it follows that
$\Sigma < .96 < S$.

Suppose $d_1 = 3$.  Then $\kappa_1 \leq \zeta_3(3) < .36$.
If $v(x_1) = 3$, then $d_i \geq 4$ and $\kappa_i \leq \zeta^3(d_i) 
< .26$, $i= 2,3$, by Lemma~\ref{LemmaX},
whence $\sum \kappa_i < S$.
Therefore $v(x_1) \geq 4$ and $\kappa_1 \leq \zeta_4(3) < .342$, so $\kappa_2
+ \kappa_3 \geq S - \kappa_1 > .6295$.  
If $d > 3$ and $d$ is odd, then $\zetao(d) < .21$.  
For all $d\geq 3$ we have $\zeta_4(d) \leq \zeta_4(3) < .342$.
It follows that $d_i$ is even whenever $d_i > 3$.  
If $d_2 > 3$, then Lemma~\ref{translation.lemma} with $\ue = (3,2,2)$
implies that
$v(x_i^2) > 1$ for some $i > 1$. 
Therefore $\kappa_2 + \kappa_3 
\leq \zeta_{4,2}(d_i) + \zeta_4(d_{5-i}) \leq \zeta_{4,2}(4) + \zeta_4(4)
< S- \kappa_1$.
This implies that $d_2 = 3$, so $d_3 > 3$.  
From Lemma~\ref{translation.lemma} with $\ue = (3,3,2)$ and $i^\ast = 3$
we have $v(x_3^2) \geq 2$.
Therefore $\sum \kappa_i \leq 2 \zeta_4(3) + \zeta_{4,2}(d_3)
\leq 2 \zeta_4(3) + \zeta_{4,2}(4) < S$, a contradiction.

\medskip
We have $d_1 = 2$ and $\kappa_1 \leq \zeta_3(2) < .5186$.
Since $\zetao(d) < .22$ for $d > 8$ it follows that $d_2 \leq 8$.
By Lemma~\ref{Scott3}, $v(x_i) \geq 5$ for $i = 2,3$.

We claim that if $i = 2$ or $3$ and $d_i > 4$, then $\kappa_i \leq .236$ and
furthermore, either $\kappa_i < .204$ or $d_i = 6$ and $v(x_i^2) \geq 3$.
Since $\zetao(d) < .2$ for $d \geq 13$ and $\zetao_5(d) < .204$ for
$d$ odd with $4 \leq d < 12$, it suffices to assume that $d_i$ is even
and $d_i \leq 12$.
We have $\kappa_i \leq \zeta_5(d_i) < .236$.
Suppose $v(x_i^2) = 1$.
By Lemma~\ref{LemmaY},
$\kappa_{5-i} \leq \zeta^1(d_{5-i})$ and $d_{5-i} \geq 11$, so 
$\kappa_{5-i} < .11$. 
It follows that $\sum \kappa_i < .97 < S$.  
Therefore $v(x_i^2) > 1$.
Suppose $v(x_i^2) = 1$.
Then $\kappa_i \leq \zeta_{5,2}(d_i) < .281$.
By Lemma~\ref{LemmaY},
$\kappa_{5-i} \leq \zeta^2(d_{5-i})$ and $d_{5-i} \geq 6$, so 
$\kappa_{5-i} < .17$. 
This also implies that $\sum \kappa_i < .97 < S$.  
Therefore $v(x_i^2) \geq 3$
and $\kappa_i \leq \zetao_{5,3}(d_i)$.
The claim follows.

It follows from the claim that if $d_2 > 4$ then $\ud = (2,6,6)$ and
$v(x_i^2) \geq 3$ for $i= 2, 3$.
By Lemma~\ref{translation.lemma} with $\ue = (2,3,3)$
we have $v(x_i^3) >1$
for some $i > 1$, so $\kappa_2 + \kappa_3 \leq  \zeta_{5,3}(6) + \zeta_{5,3,2}(6)
< .435 < S - \kappa_1$.
This shows that $d_2 \leq 4$.

Suppose $d_2 = 4$.  
Set $v = v(x_2^2)$.
If $v= 1$, then $\kappa_2 \leq \zeta_5(4) < .336$ and, as above,
$\kappa_3 < .11$.
If $v = 2$, then $\kappa_2 \leq \zeta_{5,2}(4) \leq .28$
and $\kappa_3 \leq \zeta^2(d_3) \leq \zeta^2(6) < .17$.
In either case, $\kappa_2 +\kappa_3 < .45 < S - \kappa_1$.
Therefore $v \geq 3$ and
we have $\kappa_2 \leq \zeta_{5,3}(4) < .2614$.
If $d_3 \neq 5,6,8,9,12$, then $\kappa_3 < \zetao(d_3) < .15$,
so we may assume that $d_3 \in \{5,6,8,9,12\}$.
By Lemma~\ref{24d} and the condition that $v(x_3) \geq 5$, 
$\kappa_3 \leq \zeta_{5,4,2}(d_3)$.
By inspection, this is at most $.191$ for $d_3 > 5$, so $\sum \kappa_i
< .971 < S$ in this case.  We must have $\ud = (2,4,5)$.
Thus, $A(\ud) = 2.05$ and $S = .9793$.
If $v(x_1) = 3$, then $\Sigma \leq \zeta_3(2) + \zeta^3(4) + \zeta^3(5) < .975 < S$.
Therefore $v(x_1) \geq 4$ and $\kappa_1 \leq \zeta_4(2) < .507$.
We have
$\kappa_2 \leq \zeta_{5,3}(4) < .262$ and
$\kappa_3 \leq \zeta_5(5) \leq .204$, so $\sum \kappa_i < .973 < S$,
a contradiction.
This shows that $d_2 \neq 4$, so $d_2 = 3$.

We have $S > .9781$ by Lemma~\ref{S.bounds}.
By Lemma~\ref{Scott3}, $v(x_1) \geq 4$ and $v(x_2) \geq 5$.
Since $v(x_1) + v(x_2) \geq 10$, we have
$\kappa_1 + \kappa_2 \leq \max(\zeta_4(2) + \zeta_6(3),
\zeta_5(2) + \zeta_5(3)) < .8405$.
By Lemma~\ref{23d}, $\kappa_3 \leq \zeta_{5,5,4,2}(d_3)$. 
If $d_3 > 8$, then 
$\kappa_3 < .137 < S-\kappa_1 -\kappa_2$.
Therefore $d_3 = 7$ or $8$.  

If $d_3 = 7$, then $S > .9795$.  If $n > 12$, then $v(x_1) \geq 5$, $v(x_2) \geq 7$,
and $v(x_3) \geq 7$, so $\sum \kappa_i \leq \zeta_5(2) + \zeta_7(3) + \zeta_7(7) < S$.
Therefore $n \leq 12$.
Since $\zeta_6(2) + 1/3 + 1/7 > S([2,3,7])$ we must have
$v(x_1) \leq 5$.  We have $n \geq 10$.  Therefore
$v(x_1) \leq n-v(x_1)$.
From the strong form of Scott's Theorem we have
$\max(v(x_1), n-v(x_1)) + v(x_2) + v(x_3) \geq 2n$.
Therefore $v(x_2) + v(x_3) \geq n + v(x_1) \geq 4n/3$.
Since $p=3$, we have $v(x_2) \leq 2n/3$, so $v(x_3) \geq 2n/3$.
Since $3$ has multiplicative order $6$ modulo $d_3 = 7$,
$v(x_3)$ is necessarily a multiple of $6$. 
Since $10 \leq n \leq 12$ we must have $v(x_3) = n = 12$. 
If $v(x_1) \geq 5$, then $v(x_2) \geq 7$ and
$\sum \kappa_i \geq \zeta_5(2) + \zeta_7(3) + 
\zeta_{12}(7) > S([2,3,7])$, a contradiction.
Therefore $v(x_1) = 4$ and $v(x_2) = 8$.

Suppose $d_3 = 8$. 
If $n > 10$, then $v(x_1) \geq 4$,
$v(x_2) \geq 6$,
$v(x_2^2) \geq 6$, and $v(x_2^4) \geq 3$,
so $\Sigma \leq \zeta_4(2) + \zeta_5(3) + \zeta_{6,6,1,3}(8) < S$.
Therefore $n=10$.
If $d_1 > 4$, then $d_1 = 5$, $5 \leq d_2 \leq 6$,
and $d_3 \geq 8$ by the strong form of Scott's Theorem,
so $\Sigma \leq \zeta_5(2) + \zeta_5(3) + \zeta_{8,5,1,2}(8) < S$.
Therefore $d_1 = 4$, whence $d_2 = 6$.
Since $\Sigma \leq \zeta_4(2) + \zeta_6(3) + \zeta_{6,5,1,3}(8) < S$,
we also have $v(x_3^4) = 2$.
\eop

\begin{prop}
\label{2prop}
If $p= 2$, then 
$14 \leq n \leq 21$ and one of the following is true.
\begin{enumerate}
\item $\ud = (2,3,7)$
\item $n=16$, $\ud = (2,4,5)$, $v(x_1) = 4$, $v(x_2) = 12$, and $v(x_3) = 16$.
\end{enumerate}
\end{prop}

\pf
Assume that $p=2$.
By Lemma~\ref{p.prop},
$n \geq 14$, $r \leq 8$, and
$S > (r-3) + .9589$.

\begin{step}
\label{step1}
\begin{enumerate}
\item $\zetao(2) = .75$.
\item If $d > 2$, then $\zetao(d) \leq .5$.
\item If $d > 4$, then $\zetao(d) \leq .375$.
\item If $d > 6$, then $\zetao(d) < .282$.  
\item If $d > 8$, then $\zetao(d) \leq .25$.  
\item If $d > 12$, then $\zetao(d) < .19$.
\item If $d > 14$, then $\zetao(d) \leq .15$.
\item If $d > 30$, then $\zetao(d) < .094$.  
\item If $d > 42$, then $\zetao(d) < .08$.  
\end{enumerate}
\end{step}

In view of Lemma~\ref{Lemma5},
the assertions follows immediately
from inspection of the values of $\zetao(d)$ for $d < 100$.

\begin{step}
$r < 5$.
\end{step}

If $r = 8$, then $v(x_i) \geq 2$ for at least $2$ choices of $x_i$,
so $\sum \kappa_i \leq 6 \zeta(2) + 2 \zeta_2(2) =  5.75$.
If $r = 7$, then $v(x_i) \geq 3$ for at least $2$ choices of $x_i$
since $v(x_1) + \ldots + v(x_6) \geq 14 > 6 \cdot 2$.
Therefore $\sum \kappa_i \leq 5 \zeta(2) + 2 \zeta_3(2) \leq 4.875 < S$.
This shows that $r \leq 6$.

Suppose $r = 6$.
Set $w = v(x_1) + v(x_2) + v(x_3) + v(x_4)$.
If $w \leq 6$, then Lemma~\ref{LemmaX} implies that
$d_5, d_6 > 2$ and $\kappa_i \leq \zeta^6(d_i) < .34$, $i = 5,6$,
so $\sum \kappa_i \leq 4 \zeta(2) + 2 \cdot .34 < 3.7$.
Therefore $v(x_1) + v(x_2) + v(x_3) + v(x_4) \geq 7$, and the same is
true for any other choice of $4$ distinct subscripts.  
If $v(x_i) = 1$
for $3$ values of $i$, then $v(x_j) \geq 4$ for all other values
and $\sum \kappa_i \leq 3 \zeta(2) + 3 \zeta_4(2) < 3.9$.
If $v(x_i) = 1$ for exactly $2$ values of $i$, then $v(x_j) \geq 3$
for at least $3$ values of $j$ and $\sum \kappa_i \leq 2 \zeta(2) 
+ \zeta_2(2) + 3 \zeta_3(2) < 3.9$.  
It follows that $v(x_i) = 1$
for at most $1$ choice of $i$, and $\sum \kappa_i \leq \zeta(2) + 
5\zeta_2(2) < 3.9$.
Therefore $r < 6$.

Suppose $r = 5$.  We claim that if $i, j$, and $k$ are distinct,
then $v(x_i) + v(x_j) + v(x_k) \geq 7$. 
Assume that $v(x_i) + v(x_j) + v(x_k) \leq 6$.  Then, by Lemma~\ref{LemmaX},
$d_l > 2$ for
$l \neq i,j,k$ and $\kappa_l \leq \zeta^6(d_l)$.
If $d_l > 6$, then 
$\kappa_l < .3$ by Step~\ref{step1}.  
If $3 \leq d_l \leq 6$, then $\zeta^6(d_l) < .34$ by inspection.
This implies that $\sum \kappa_i < 3 \zeta(2) + 
2 \cdot .34 = 2.93 < S$, and the claim follows.

We claim further that if $v(x_i) + v(x_j) \leq 4$ for distinct $i,j$,
then $d_k = 2$ for all $k \neq i,j$.  
For the purpose of establishing this claim we remove the running
assumption on the ordering of $x_i$ for the balance of this paragraph
and show that if $v(x_1) + v(x_2) \leq 4$ then $d_k = 2$ for $k > 2$.
If $v(x_1) + v(x_2) = 2$, then
$v(x_k) \geq 5$ for $k > 2$ 
by the previous paragraph,
and $\sum \kappa_i \leq 2 \zeta(2) + \sum_{k > 2} \zeta_5(d_k)$.
Since $\zeta_5(2) < .52$ and $\min(\zeta_5(d),\zetao(d)) < .4$
for $d > 2$, we have either $\sum \kappa_i \leq
2 \zeta(2)  + 2 \zeta_5(2) + .4 < 2.94$ or $d_k =2$
for all $k > 2$.
If $v(x_1) + v(x_2) = 3$, then $\kappa_1 + \kappa_2 \leq \zeta(2)
+ \zeta_2(2) = 1.4$ and $\sum \kappa_i \leq \kappa_1 + \kappa_2 + 
\sum_{k > 2} \zeta_4(d_k)$.  
Since $\zeta_4(2) < .54$ and $\min(\zeta_4(d),
\zetao(d)) < .41$ when $d > 2$, either $\sum \kappa_i < 2.9$
or $d_k = 2$ for all $k > 2$. 
Finally, if $v(x_1) + v(x_2) = 4$, then $\kappa_1 + \kappa_2 \leq 
\zeta(2) + \zeta_3(2) < 1.32$.  
Considering that $\zeta_3(2) < .57$ and $\min (\zeta_3(d), \zetao(d))
< .44$ for $d > 2$, either $\sum \kappa_i < 2.9$
or $d_k = 2$ for all $k > 2$.
Since $\sum \kappa_i \geq S$ we conclude in
every case that $d_k = 2$ for all $k > 2$.
This completes the argument that if $v(x_i) + v(x_j) \leq 4$
for some $i \neq j$, then $d_k = 2$ whenever $k \neq i, j$.

Reverting to the ordering of $x_i$, so that $d_5$ is the largest
value of $d_i$,
the previous paragraph implies that
if $d_5 > 2$, then  $v(x_i) + v(x_j) \geq 5$ for every pair of distinct 
$i, j < 5$.
In that case,
$\sum_{i < 5} \kappa_i \leq \max(\zeta(2) + 3 \zeta_4(2),
\zeta_2(2) + 3 \zeta_3(2)) < 2.4$, and $\kappa_5 \leq \zetao(d_5) \leq .5$,
whence $\sum \kappa_i < S$.
We conclude that $d_i = 2$ for all $i$.
From Lemma~\ref{translation.lemma} with $\ue = (2,2,2,1,1)$
it follows that $v(x_i) + v(x_j) \geq 7$
whenever $i \neq j$, whence  
$\kappa_i \leq \max(\{\zeta_a(2) + 4 \zeta_{7-a}(2)
\ : \ a = 1,2,3\}) < 2.8 < S$.
This completes the argument that $r \neq 5$.

\begin{step}
\label{step3}
$r = 3$.
\end{step}

Suppose $r = 4$.  Since $v(x) + v(x') + v(x'') \geq 14$ for every set
of $3$ generators $\{x,x',x''\}$ it follows that $v(x) \geq 5$
for at least two of the four generators, so $\kappa_i \leq \zeta_5(d_i)$
for at least two values of $i$.  

We claim that $A(\ud) \leq 3$. 
Suppose $A(\ud) > 3$.
Then $\sum 1/d_i < 1$. 
The ordering assumption on $d_i$ implies that $d_2 > 2$ and $d_4 > 4$, 
so $\kappa_i \leq .5$ for $i > 1$ and $\kappa_4 \leq .375$.
If $d_1 > 2$, then $\sum \kappa_i \leq 1.875$, which is not the
case, so $d_1 = 2$.  It follows that $d_3 > 4$, since 
otherwise $1/d_1 + 1/d_2 + 1/d_3 \geq 1$.
This implies that
$\kappa_1 + \kappa_2 + \kappa_3
\leq 1.625$.  Therefore $\kappa_4 > .3$, so $d_4 \leq 6$ and
$A(\ud) \leq A(2,4,6,6) < 3$, a contradiction.
This establishes the claim, and we conclude that $S \geq 1.9698$.

Set $w = v(x_1) + v(x_2)$. 
Then, by Lemma~\ref{LemmaX}, $\kappa_3 \leq \zeta^w(d_3)$,
$\kappa_4 \leq \zeta^w(d_4)$, and $d_3 \geq 14/w$.
If $w \leq 3$, then $\kappa_1 + \kappa_2 \leq 1.5$, $d_i \geq 5$,
and $\kappa_i < \zeta^3(d_i) < .201$ for $i > 2$.
If $w = 4$, then $\kappa_1 + \kappa_2 \leq \zeta(2) + \zeta_3(2) < 1.32$,
$d_i \geq 4$, and $\kappa_i < \zeta^4(d_i) < .26$ for $i > 2$.
If $w = 5$, then $\kappa_1 + \kappa_2 \leq \zeta(2) + \zeta_4(2) < 1.282$,
$d_i \geq 3$, and $\kappa_i < \zeta^5(d_i) < .335$ for $i > 2$.
If $w = 6$, then $\kappa_1 + \kappa_2 \leq \zeta(2) + \zeta_5(2) < 1.266$,
$d_i \geq 3$, and $\kappa_i < \zeta^6(d_i) < .336$ for $i > 2$.
In each case, $\sum \kappa_i \leq 1.96 < S$.
This implies that $v(x_1) + v(x_2) \geq 7$.
More generally, 
$v(x_i) + v(x_j) \geq 7$ whenever $i \neq j$.

Suppose $d_3 > 2$ and
set $v = v(x_1)$.
We claim that $v = 1$.
If $v = 2$, then $\kappa_i \leq \zeta_5(d_i)$ for all $i > 1$.
Since $\zetao(d) < \zeta_5(4) < .5 $ when $d > 4$ and
$\zetao_k(3) \leq \zetao_k(4)$ for all $k$ it follows that
$\Sigma \leq \zeta_2(2) + \zeta_5(2) + 2 \zetao_5(4) < 1.93$.
Similarly, if $v = 3$, then $\Sigma \leq \zeta_3(2) + \zeta_4(2) + 2 \zetao_4(4) < 1.91$.
If $v = 4$, then $\Sigma \leq 2\zeta_4(2) + \zetao_3(4) + \zetao_4(4) < 1.91$.
If $v = 5$, then $\Sigma \leq 2\zeta_5(2) + \zetao_2(4) + \zetao_5(4) < 1.93$.
If $v \geq 6$, then $\kappa_1 + \kappa_2 \leq 2 \zeta_6(2) < 1.02$
and $\kappa_3 + \kappa_4 \leq \max_{t = 1,2,3}(\zetao_t(4) + \zetao_{7-t}(4)) < .9$.
In all cases, $\Sigma < S$.

Therefore $v(x_1) = 1$ and
$v(x_i) \geq 6$ for $i > 1$. 
We have $\kappa_1 + \kappa_2 \leq \zeta(2) + \zeta_6(2) < 1.26$.
Also, $\kappa_i \leq \zetao_6(d_i)$ when $i > 2$.
If $d > 4$, then $\zetao_6(d) < .34$.
Therefore $d_3 \leq 4$ and $\kappa_3 < .39$.
From Lemma~\ref{translation.lemma} with $\ue = (1,2,d_3,2)$
we have
$2d_3 + d_3 v(x_4^2) \geq 28$, whence $v(x_4^2) \geq 28/d_3 -2 \geq 5$.
Therefore $\kappa_4 \leq \zetao_{6,5}(d_4)$.
If $d_4 \geq 4$, then $\kappa_4 < .3$ and
$\Sigma < 1.95$.  Therefore $d_4 =  3$, whence $d_3 = 3$, and
$\kappa_i \leq \zeta_6(3) < .35$ for $i = 3$ or $4$.  Once again,
$\Sigma < S$.
This shows that $d_3 = 2$.

We have $A(\ud) < 2.5$, so $S > 1.9748$.
As before, set $v = v(x_1)$.
From Lemma~\ref{LemmaX4}, $\kappa_4 \leq \zeta^{2v}(d_4)$ and $d_4 \geq 7/v$.
If $v = 1$, then $\kappa_1 + \kappa_2 + \kappa_3 \leq \zeta(2) + 2\zeta_6(2)
< 1.766$ and $\kappa_4 \leq \zeta^2(d_4) < .144$ because $d_4 \geq 7$.
If $v = 2$, then $\kappa_1 + \kappa_2 + \kappa_3 \leq \zeta_2(2) + 2\zeta_5(2)
< 1.657$ and $\kappa_4 \leq \zeta^4(d_4) < .255$ because $d_4 \geq 4$.
If $v = 3$, then $\kappa_1 + \kappa_2 + \kappa_3 \leq \zeta_3(2) + 2\zeta_4(2)
= 1.625$ and $\kappa_4 \leq \zeta^6(d_4) < .336$ because $d_4 \geq 3$.
This shows that $\Sigma < S$ when $v \leq 3$.  Therefore $v \geq 4$
and $\kappa_1 + \kappa_2 + \kappa_3 \leq 3 \zeta_4(2) < 1.594$.
From Lemma~\ref{translation.lemma} with $\ue = (2,2,2,1)$
we have $v(x_4) \geq 7$, so $\kappa_4 \leq 
\zeta_7(d_4) < .379$ because $d_4 > 2$.
In this case as well, $\sum \kappa_i < S$.  This shows that $r < 4$.

\begin{step}
\label{leq4}
$v(x_i) \geq 4$ for all $i$.
\end{step}

Since $A(\ud) < r = 3$, $S > .9698$.
Set $v = v(x_1)$.  
We apply Lemma~\ref{LemmaX} once again to bound $v$ from below.
If $v = 1, 2$, or $3$, then $\kappa_1 \leq \zeta_v(2)
\leq .75, \ .625, \ .563$, respectively.
For $i > 1$, $\kappa_i \leq \zeta^v(d_i)$ where $d_i \geq 14, 7, 5$
in the respective cases.
Using Step~\ref{step1} and inspection, we have
$\kappa_i < .08, .15, .201$ in the respective cases.
It follows that $\sum \kappa_i < S$ whenever $v(x_1) < 4$.
Therefore $v(x_1) \geq 4$.
More generally, since the argument that established this
does not use the ordering assumption
on $x_i$, it follows that $v(x_i) \geq 4$ for all $i$.

\begin{step}
\label{step5}
$d_1 = 2$.
\end{step}

Assume that $d_1 > 2$.  
It follows from Step~\ref{step1} that $d_1 \leq 6$,
so $A(\ud) < 2.84$ and $S > .9714$.
Since $v(x_1) \geq 4$, we have  $\kappa_1 \leq \zetao_4(d_1)$. 
It follows from Step~\ref{step1} and inspection that $\kappa_i
< .41$ for all $i$.

If $v(x_1) = 4$, then $d_i \geq 4$, $i = 2, 3$ and
$\kappa_i \leq \zeta^4(d_i) < .255$, which implies that $\sum \kappa_i < S$.
Therefore $v(x_1) \geq 5$, and, similarly, $v(x_i) \geq 5$ for all $i$.
Thus $\kappa_i \leq \zetao_5(d_i)$ for all $i$.  
In particular, $\kappa_i < .3907$ for all $i$.

Suppose $d_1 > 4$.  Since $\zetao_5(d) < .27$ when $d> 4$, $d \neq 6$
and $\zetao_5(6) < .35$, it follows that 
$d_i = 6$ for all $i$.
Lemma~\ref{zeta5*}
implies that $v(x_i^2) \geq 3$ for at least two choices of $i$,
so $\sum \kappa_i \leq \zetao_5(6) + 2 \zetao_{5,3}(6)< .95$.
Therefore $d_1 \leq 4$.

\smallskip
Suppose $d_1 = 4$.  
Then $d_2 \leq 6$ since otherwise $\kappa_i \leq \zetao_5(d_i) < .27$
for $i = 2, 3$ and $\Sigma < \zetao_5(4) + 2 \cdot .27 < S$.
It follows from Step~\ref{step1} that $d_3 \leq 12$ as otherwise
$\Sigma < S$.
This implies that
$A(\ud) \leq A(4,6,12) = 2.5$, so $S \geq .9748$.
Also, Lemma~\ref{zeta5*} implies that $v(x_1^2) + v(x_2^2) \geq 3$.
We claim that $d_3 \leq 8$.  
If $d_2 = 5$ or $6$, then $\kappa_1 + \kappa_2 < \zetao_5(4) + \zetao_5(6) < .7345$,
so $\zetao_5(d_3) \geq \kappa_3 > .24$.
It follows from inspection that $d_3 \leq 8$ in this case.
If $d_2 = 4$, then
$\kappa_1 + \kappa_2 \leq \zetao_5(4) + \zetao_{5,2}(4) < .7188$,
so $\zetao_5(d_3) > .25$ and $d_3 \leq 8$ in this case as well.

From Lemma~\ref{zeta5*}.\ref{part1} we have $v(x_1^2) + v(x_2^2) \geq 4$ and
$v(x_1^2) + v(x_3^2) \geq 5$.
If $v(x_1^2) = 1$, then $v(x_2^2) \geq 3$ and $v(x_3^2) \geq 4$. 
so $\kappa_2 \leq \zetao_{5,3}(d_2) < .3021$,
$\kappa_3 \leq \zetao_{5,4}(d_3) < .2813$,
and $\sum \kappa_i < .974 < S$.
If $v(x_1^2) = 2$, then $\kappa_1 \leq \zetao_{5,2}(4) < .3282$,
$\kappa_2 \leq \zetao_{5,2}(d_2) < .3438$, and $\kappa_3 \leq \zetao_{5,3}(d_3)
< .3021$,
whence $\sum \kappa_1 < .9741 < S$.
We conclude that $v(x_1^2) \geq 3$, so that $\kappa_1 \leq \zeta_{5,3}(4) < .3$.
Without loss, if $d_i = 4$, $i = 2, 3$, then $\kappa_i \leq
\zetao_{5,3}(4) < .3$.  
If $d_i > 4$ for some $i$, then $\kappa_i \leq \zetao_5(d_i) < .35$.
Since $v(x_2^2) + v(x_3^2) \geq 7$ by Lemma~\ref{zeta5*}, we have either
$v(x_2^2) \geq 4$ or $v(x_3^2) \geq 4$, whence $\kappa_i \leq \zetao_{5,4}(d_i) < .3$
for some $i > 1$.
It follows that $\Sigma < .95 < S$, so we conclude that $d_1 \neq 4$.

\medskip
We may therefore suppose $d_1 = 3$, so that $\kappa_1 \leq \zetao_5(3) \leq .3542$.
By Lemma~\ref{Scott3}, $\kappa_i \leq \zetao_5(d_i)$ for $i = 2,3$.
As in the argument when $d_1 = 4$, it follows that $d_2 \leq 6$.
By Lemma~\ref{translation.lemma}
with $\ue = (1,1,1)$, we have $v(x_i) \geq 7$ for two choices of $i$.
If $d_2 = 6$, then $\kappa_1 + \kappa_2 \leq \max(\zetao_5(3)+\zetao_7(6),
\zetao_7(3)+\zetao_5(6)) < .6902$.  
It follows from inspection of $\zetao_5$ values that $d_3 = 6$.
Since $v(x_2^2) + v(x_3^2) \geq 10$ by Lemma~\ref{zeta5*}.\ref{part1}
we have $\sum \kappa_i \leq \zetao_5(3) + \zetao_5(6) + \zetao_{5,5}(6)
< .3542 + .3438 + .2709 < .97 < S$.

If $d_2 = 5$, then $\kappa_2 \leq .225$ and $\kappa_3 \leq \zetao_5(d_3) < .344$,
so $\Sigma < S$.

If $d_2 = 4$, then $\kappa_1 + \kappa_2 \leq \max(\zetao_5(3) + \zetao_7(4),
\zeta_7(3) + \zeta_5(4) ) < .7332$.
By Lemma~\ref{zeta5*}.\ref{part3}, $\kappa_3 \leq \zetao_{5,3}(d)$.
It follows that $\zetao_{5,3}(d) > .24$, so $d_3 = 4$ or $6$ by inspection.
If $d_3 = 4$, then $\kappa_2 \leq \zetao_{5,3}(4) < .3$ by the same result,
and $\sum \kappa_i < \zetao_5(3) + 2 \zetao_{5,3}(4) < S$.
Therefore $d_3 = 6$.
If $v(x_2^2) \geq 3$, then 
$\sum \kappa_i \leq \zetao_5(3)
+ \zetao_{5,3}(4) + \zetao_{5,3}(6) < S$.
If $v(x_2^2) = 2$, then $v(x_3^2) \geq 8$ by Lemma~\ref{zeta5*}
and $\sum \kappa_i < \zetao_5(3)
+ \zetao_{5,2}(4) + \zetao_{5,5}(6) < S$, so
we may assume that $v(x_2^2) = 1$.
From Lemma~\ref{translation.lemma} with $\ue = (3,2,3)$
we have $4v(x_3^3) + 6 \geq 28$
whence $v(x_3^3) \geq 6$ and $\kappa_3 < \zetao_{5,5,6}(6)
< .2 < S - \kappa_1 - \kappa_2$.  This shows that $d_2 \neq 4$.

If $d_2 = 3$ then $\kappa_2 \leq \zetao_7(3) \leq .3386$ and, by
Lemma~\ref{zeta5*}.\ref{part2},
 $\kappa_3 \leq \zetao_{5,5}(d_3) \leq .2735$, so $\Sigma 
\leq .97 < S$.

\begin{step}
\label{step6}
$d_2 \leq 4$
\end{step}

By Lemma~\ref{S.bounds} and the previous step, $S \geq .9748$.
Assume that $d_2 > 4$.  By Step~\ref{leq4}, $v(x_1) \geq 4$.
If $v(x_1) = 4$, then $\kappa_1 \leq \zeta_4(2) < .532$,
and $\kappa_i \leq \zetao_{10,6,2}(d_i)$ by Lemma~\ref{LemmaX}.
By inspection, $\kappa_i < .22$ for $i \geq 2$.
This implies that $\Sigma < S$.  We conclude that $v(x_1) > 4$.

We have $\kappa_1 \leq \zeta_5(2) < .5157$.  If $d_i > 8$ and $d_i \neq 12$,
then $\kappa_i \leq \zetao(d_i) < .2$.  
If $d_i = 12$, then $\kappa_i \leq \zetao_7(12) < .232$.  
It follows that  either $d_2 \leq 8$ or $d_2 = d_3 = 12$.
In the latter case, Lemma~\ref{translation.lemma} with $\ue = (2,3,3)$
shows that $v(x_2^3) + v(x_3^3) \geq 7$,
so $v(x_i^3) \geq 4$ for some $i > 1$, and 
$\kappa_i \leq \zetao_{7,1,4}(12) < .21$.
This implies that $\Sigma < S$.  
We conclude that $d_2 \leq 8$.
From Lemma~\ref{translation.lemma} with $\ue = (2,d_2,2)$
we have $v(x_3^2) \geq 2 \cdot 14/d_2 > 3$.
Consequently, $\kappa_3 \leq \zetao_{7,4}(d_3)$.
Suppose $d_2 = 8$.  
Then $\kappa_2 \leq \zetao_7(8) < .254$. If $d_2 > 12$, then
$\kappa_3 < .19$ by Step~\ref{step1} and
$\Sigma < S$, so $d_2 \leq 12$.
By Lemma~\ref{translation.lemma} with $\ue = (2,2,12)$, 
$v(x_2^2) > 2$,  
so $\kappa_2 \leq \zetao_{7,3}(8) < .223$.
Since $\zetao_{7,4}(12) < .222$, we conclude that $\kappa_2 + \kappa_3
< .446 < S - \kappa_1$, a contradiction. 
Therefore $d_2 < 8$.
Since $\zetao_7(7) < .15$, it is evident that $d_2 \neq 7$.

Suppose $d_2 = 6$.
Then $A(\ud) < 2.34$ and $S > .9764$, so $\kappa_2 + \kappa_3 \geq S- \kappa_1
> .4607$.
Set $w =v(x_2^2)$.
Then $w$ is necessarily even because $x_2^2$ has order $3$.
If $w = 2$, then $\kappa_2 \leq \zetao_7(6) < .336$. 
By Lemma~\ref{LemmaY}, $d_{3} \geq 14$ and $\kappa_3 \leq \zeta^1(d_3)$.
By Step~\ref{step1} and inspection of the values of $\zeta^1(d)$
for $14 \leq d \leq 30$ we have $\kappa_3 < .08$, 
so $\Sigma < S$ in this case.
If $w= 4$, then $\kappa_2 \leq \zetao_{7,4}(6) < .2735$.
By Lemma~\ref{LemmaY}, $d_3 \geq 7$
and $\kappa_3 \leq \zeta^2(d_3)$.
Observing that $\zeta^2(d) < .1431$ for $7 \leq d \leq 28$,
we conclude from Step~\ref{step1} that $\Sigma < S$
in this case as well.
If $w = 6$, then $\kappa_2 \leq \zetao_{7,6}(6) < .2579$.
We have $d_3 \geq d_2 = 6$, and, by Lemma~\ref{LemmaY},
$\kappa_3 \leq \zeta^3(d_3)$.
Since $\zeta^3(d) < .18$ for $6 \leq d \leq 12$ we conclude
from Step~\ref{step1} that $\kappa_3 < .18$, whence, once again,
$\Sigma < S$.
It follows that $w \geq 8$, so $\kappa_2 \leq \zetao_{8,8}(6) < .2527$.
From Lemma~\ref{translation.lemma} with $\ud = (2,6,2)$
we have $v(x_3^2) \geq 5$, so
$\kappa_3 \leq \zetao_{7,5}(d_3)$.
If $d_3 > 6$ and $d_3 \neq 12$, then $\zetao_{7,5}(d_3) < .2$ and
$\Sigma < S$.  Therefore either $d_3 = 6$ or $d_3 = 12$.  
Recall that, by Lemma~\ref{translation.lemma} with $\ue = (2,3,3)$,
$v(x_2^3) + v(x_3^3) \geq 7$.
If $v(x_2^3) = 1$, then $\kappa_3 \leq \zetao_{7,5,6}(d_3) < .2$,
and $\Sigma < S$.
Therefore $v(x_2^3) \geq 2$, so $\kappa_2 < \zetao_{8,8,2}(6)
< .211$.  If $d_3 = 6 = d_2$, then we may assume that $\kappa_3 \leq \kappa_2$,
whence $\kappa_2 + \kappa_3 < .43$.
If $d_3 = 12$, then $\kappa_3 \leq \zetao_{7,5}(12) < .217$
and $\kappa_2 + \kappa_3 < .43$.
In either case, $\Sigma < S$.   Therefore $d_2 \neq 6$.

Suppose $d_2 = 5$.
Then $\kappa_2 < .2063$, so $\kappa_3 \geq S - \kappa_1 - \kappa_2 > .25$.
We have $\kappa_3 < \zetao_{7,6}(d_3)$ by
Lemma~\ref{translation.lemma} with $\ue = (2,5,2)$.
It follows from Step~\ref{step1} and inspection that $d_3 = 6$.
From Lemma~\ref{translation.lemma} with $\ue = (2,5,3)$
we have $v(x_3^3) \geq 2$,
so $\kappa_3 < \zetao_{7,6,2}(d_3) < .22$, a contradiction.

\begin{step}
\label{step7}
If $d_2 = 4$, then 
$n = 16$, 
$\ud = (2,4,5)$,
$v(x_1) = 4$, $v(x_2) = 12$, and $v(x_3) = 16$.
\end{step}

Suppose $d_2 = 4$.
Then $A < 2.25$ and $S > .9773$.
Also, $\kappa_2 \leq \zetao_7(4) < .379$.

Assume that $v(x_1) = 4$.
then $\kappa_1 \leq \zeta_4(2) < .532$.
By Lemma~\ref{LemmaX}, $d_3 > 3$ and $\kappa_2 \leq \zeta^4(d_3)$,
so $\kappa_2 < .255$ by inspection and Step~\ref{step1}.
From Lemma~\ref{translation.lemma} with $\ue = (1,4,4)$
we have $v(x_3^4) \geq 2n - 4 \cdot 4 \geq 12$, so
$v(x_3) \geq 12$ and $v(x_3^2) \geq 12$ as well.
By Lemma~\ref{24d} we have $v(x_3^3) \geq 4$.
Therefore $\kappa_3 \leq \zetao_{12,12,4,12}(d_3)$.
If $d_3 > 5$ then $\kappa_3 < .178$ by Step~\ref{step1}
and inspection.  Therefore $d_3 = 5$.
We have $n \leq d_2 v(x_1) \leq 16$ and $v(x_i) \geq n - 4 \geq 10$, $i = 2,3$.
Since $2$ has multiplicative order $4$ modulo $5$, we also have
$4  |  v(x_3)$, so $v(x_3) = 12$ or $16$. 
If $v(x_3) = 12$, then $v(x_2) \geq 2n - v(x_1) - v(x_3) = 16$.
However, $v(x_2) \leq 3n/4$ because $x_2$ is an element of order $4$
acting in characteristic $2$. 
These inequalities are not compatible
with the condition $n \leq 16$.  We conclude that $v(x_3) = 16$, $n= 16$,
and $v(x_2) = 12$.

We may therefore assume that $v(x_1) > 4$.
Then $\kappa_1 \leq \zeta_5(2) < .5157$.
Set $w = v(x_2^2)$.
Assume that $w \leq 2$. 
Then, by Lemma~\ref{LemmaY},
$d_3 \geq 14$, 
and $\kappa_3 \leq \zeta^1(d_3)$. 
So $\kappa_3 < .08 < S - \kappa_1 - \kappa_2$.
Therefore $w > 2$.
If $w = 3$ or $4$, then $\kappa_2 \leq \zetao_{7,3}(4) < .2852$,
$d_3 \geq 7$, and $\kappa_3 \leq \zeta^2(d_3)$, so $\kappa_3 < .144$
by inspection and Step~\ref{step1}.
Once again, $\Sigma < S$.
If $w = 5$ or $6$, then $\kappa_2 \leq \zetao_{7,5}(4) < .2618$,
and $\kappa_3 \leq \zeta^3(d_3)$. 
If $d_3 \geq 6$, then $\kappa_3 <.19$ and $\Sigma < S$, so $d_3 = 5$.
By Lemma~\ref{LemmaY}, $w = 6$.
Thus, $\kappa_2 \leq
\zetao_{7,6}(4) < .2579$ and $\kappa_3 \leq \zeta^3(5) < .2004$,
so $\Sigma < S$.
We conclude that $w = v(x_2^2) \geq 7$, 
so $\kappa_2 \leq \zetao_{7,7}(4) < .2559$.
By Lemma~\ref{24d}, $\kappa_3 \leq \zetao_{7,7,4}(d_3)$.
If $d_3 > 5$, then $\kappa_3 < .2$ by Step~\ref{step1} and inspection,
so $\Sigma < S$.
If $d_3 = 5$, then $S = .9793$, and 
$\kappa_3 \leq .2063$, so once again $\Sigma < S$.  
This completes the argument that $d_3 \neq 4$.

\begin{step}
\label{step8}
If $d_2 = 3$ then $\ud = (2,3,7)$.
\end{step}

It suffices to assume that $d_2  = 3$ and $d_3 > 7$.
We have $A(\ud) < 2.17$ and $S > .9781$.
Also, $v(x_2)$ is even because $x_2$ is an element of order $3$ acting
over $\F_2$.
In particular, $v(x_2) \geq 8$ and $\kappa_2 < .33595$. 
We have $\kappa_3 \leq \zetao_{10,10,7,5,3}(d_3)$ by Lemma~\ref{23d}.
By inspection, $\kappa_3 \leq .132$.
If $v(x_1) \geq 6$, then $\kappa_1 < .50782$ and $\Sigma < S$,
so $v(x_1) = 5$ by Lemma~\ref{Scott3}.
We have $\kappa_1 \leq \zeta_5(2) < .5157$.

It follows that $v(x_2) \geq n-5 \geq 9$, whence $v(x_2) \geq 10$, and
$\kappa_2 \leq \zeta_{10}(3) < .334$. 
We have $\kappa_1 + \kappa_2 < .8497$.

By inspection, if $d > 7$ and $d \neq 8$ or $12$, then $\zetao_{10,10,7,5,3}(d) < .114$.
It follows that $d_3 = 8$ or $12$.
If $d_3 = 8$, then $A(\ud) < 2.05$, so $S > .9793$
and $\Sigma \leq \zeta_5(2) + \zeta_{10}(3) + \zetao_{10,10,7,5}(8) < .9793 < S$.
We conclude that $d_3 = 12$, whence $A(\ud) < 2.09$ and
$S > .9789$.
Since $x_3^4$ has order $3$, $v(x_3^4)$ must be even, and
$v(x_3^4) \geq 6$.
If $v(x_2) = 10$, then $v(x_3) \geq 2n-v(x_1) - v(x_2) \geq 13$,
so $\kappa_3 \leq \zetao_{13,10,7,6}(12)$,
and $\sum \kappa_i \leq \zeta_5(2) + \zeta_{10}(3) + \zetao_{13,10,7,6}(12)$.
If $v(x_2) > 10$, then $v(x_2) \geq 12$ and
$\sum \kappa_i \leq \zeta_5(2) + \zeta_{10}(3) + \zetao_{10,10,7,6}(12)$.
In either case, $\Sigma < S$, a contradiction.

\begin{step}
Conclusion
\end{step}

By Steps~\ref{step3},
\ref{step5},
\ref{step6},
\ref{step7}, and
\ref{step8},
it suffices to show that 
if $\ud = (2,3,7)$ then $14 \leq n \leq 21$.

Assume that $\ud = (2,3,7)$.  
By Lemma~\ref{n.bounds}, $n \geq 14$.
If $n > 21$, then $v(x_1) \geq 8$, $v(x_2) \geq 11$, and $v(x_3) \geq 11$,
so $\sum \kappa_i \leq \zeta_8(2) + \zeta_{11}(3) + \zeta_{11}(7) < .9795 < S$.
\eop

Theorem~\ref{basic.result} now follows from 
Propositions~\ref{large.p}--\ref{2prop}.

Note that for $p = 2,3$, or $5$, 
further information about values of $v(y)$ for certain  
elements $y$ is recorded in Propositions~\ref{p=5.prop}, 
\ref{p=3.prop}, and~\ref{2prop}.

\section{Proof of Theorem~\ref{grassmann.exceptions}}

Retaining the notation of \ref{notation.section},
assume that $\Omega$ is a primitive point action for $G$ with $\order{\Omega} \geq 10^4$
and that $x \in G$.

\subsection{Linear and Symplectic Groups}
\begin{prop}
\label{linear}
If $\Omega$ consists of all points in the $L$ action or $Sp$ action, then
$f(x) - q^{-v(x)} < 1/100$.
\end{prop}

\pf
We have $N = (q^n-1)/(q-1)$, so $q^{n-1} < N < 2q^{n-1} \leq q^n$.

Suppose $x$ is a linear transformation.
Then the fixed points of $x$ are contained in the union of its eigenspaces,
the largest of which has dimension $n-v$.
We claim $f(x)  - q^{-v(x)} < q^{-n/2} < 1/100$.
It suffices to establish the first inequality.

If $v \leq n/2$, then the fixed points of $x$ lying outside
the largest eigenspace are contained in a space of dimension $v$.
This implies that $f(x) \leq \displaystyle
\frac{q^{n-v} -1}{q-1} + \frac{q^v-1}{q-1}$, so
$$ \begin{array}{rcl}
\displaystyle
\frac{F(x)}{N} - q^{-v} & \leq & \displaystyle \frac{q^{n-v} - 1}{q^n - 1} + 
\frac{q^v-1}{q^n-1} - q^{-v} \\
& < & q^{-(n-v)} \leq q^{-n/2}.
\end{array}$$

If $v = \frac{n+1}{2}$, then the fixed points of $x$ lying outside
the largest eigenspace are contained in the union of two nontrivial
spaces having total dimension $n-v=(n+1)/2$.  
For fixed $m$, the largest value of $q^{a} + q^{m-a}$ for $a$
in $\{1,2,\ldots,m-1\}$ is $q^{m-1} + q$.
Therefore
$\displaystyle F(x) \leq \frac{q^{n-v}-1}{q^n-1} +  \frac{q^{(n-1)/2}-1}{q-1} + 1$, so
$$\frac{F(x)}N - q^{-v} < \frac{q^{(n-1)/2} - 1}{q^n -1} + \frac{q-1}{q^n-1} 
< q^{-n/2}.$$

If $v \geq n/2 + 1$, then $x$ has at most $q-1$ eigenspaces, each of which 
has dimension at most $n/2 - 1$, so $\displaystyle F(x) \leq (q-1) 
\frac{q^{n/2-1}-1}{q-1}$ and
$$\frac{F(x)}{N} \leq (q-1)\left(\frac{q^{n/2-1} - 1}{q^n-1}\right) < q^{-n/2}.$$
This completes the analysis for $x$ a linear transformation.

\bigskip
Now suppose $x$ is not a linear transformation.  
Then $x$ induces a field automorphism because graph automorphisms do not 
act on $\Omega$.  
Let $d$ be the order of $x$ modulo InnDiag.
Then $\displaystyle F(x) \leq \frac{q^{n/d}-1}{q^{1/d} -1}$,  
so $f(x) > .01$ implies that
$$q^{n(d-1)/d} < \frac{q^n-1}{q^{n/d} -1} < 100 \frac{q-1}{q^{1/d} - 1} =
100 q^{(d-1)/d}\left( \frac{1- q^{-1}}{1 - q^{-1/d}}\right).$$
Since $q^{-1/d} \leq 1/2$, we have
$q^{n(d-1)/d} < 200 q^{(d-1)/d}(1-q^{-1}) < 200q^{(d-1)/d}$.  
It follows that $q^{(n-1)(d-1)/d} < 200$, so $200^{d/(d-1)} > q^{n-1}$.

By the first line of this argument,
$2q^{n-1} > N > 10000$. 
Therefore $q^{n-1} > 5000 > 200^{3/2}$,
whence $\frac{d}{d-1} > \frac32$, and $d = 2$.

If $x$ is not a standard field automorphism, then $F(x) \leq 
\displaystyle \frac{q^{n/2-1}-1}{q^{1/2}-1} + 1$,
so
$$
\begin{array}{rcl}
.01 < f(x) & \leq & 
(q^{1/2}+1)\left(\frac{q^{n/2-1}-1}{q^n-1}\right) + \frac1N \\
& < & \frac32q^{1/2}\cdot q^{-(n/2 + 1)} + .0001
\end{array}$$
This implies that $q^{n+1} < \displaystyle \left(\frac1{.0066}\right)^2 < 160^2$.

On the other hand, we have $F(x) > .01 N > 100$,  so 
$\displaystyle \frac{q^{n/2-1} - 1}{q^{1/2} - 1} + 1 > 100$.
It follows from this that $q^{n-2} > 99^2$, whence $q^3 < (160/99)^2$, 
which is impossible.
Therefore $x$ must be a standard field automorphism.

We have $f(x) = \frac{q^{1/2} + 1}{q^{n/2} + 1}$ and $v_q(x) = n/2$.  
If $f(x) - q^{-v_q(x)} > .01$,
then $q^{-(n-1)/2} > .01$, whence $q^{n-1} < 10000$.  
On the other hand, $q^{n-1} \cdot \displaystyle
\frac{q}{q-1} > \frac{q^n-1}{q-1} = N > 10000$.
That is,
$$q^{n-1} < 10000 < \frac{q^n}{q-1}.$$
Since $n > 2$, the first inequality implies that $q < 100$.
Since $q$ is both a perfect square and a prime power, it follow 
easily by inspection
that these two inequalities cannot both hold.
\eop

\begin{prop}
\label{Sphypprop}
If $\Omega$ consists of hyperplanes of type $\delta$ in the $Sp$ action, then
$f(x) < q^{-v(x)} + 1/100$.
\end{prop}

\pf
We have $N = \frac12(q^n + \delta q^{n/2})$.
Since $q^n$ is an even power of $2$ and $2^{14}+2^7 < 20000$, 
we have $q^n \geq 2^{16}$.

If $x$ is a field automorphism, then $F(x) \leq q^{n/2}$ in either action,
so $f(x) \leq 2(q^{n/2}-1)^{-1} < .01$.

If $x$ is in InnDiag, then $F(x) \leq \frac12(q^{n-v}+ q^{n/2})$,
so $F(x) - q^{-v(x)}N < q^{n/2}$, and $f(x) - q^{-v} < .01$,
as before.
\eop

\subsection{Actions of Unitary and Orthogonal groups}
We record here properties of orthogonal and unitary actions 
that will be used in the analysis.

\newcommand\AAA{P}
\newcommand\BBB{S}
\newcommand\cB{{\cal B}}
\begin{fact}
\label{point.count.lemma}
Let $W$ be an orthogonal or hermitian space of dimension $m$ over
$\F_q$, and let $\pi(W)$ be the number of points of a given type 
in $W$.  If $\rad W$, the totally singular radical of $W$, 
has dimension $r$, then 
$$\AAA(m) - \BBB(m+r) \leq \pi(W) \leq \AAA(m) + \BBB(m+r)$$
where $\AAA(m)$ and $\BBB(m)$ are as given below.  
$$ 
\begin{array}{ccc}
{\rm{type}} & \AAA(m) & \BBB(m) \\
\\
U,\bs & \displaystyle \frac{q^{m-1/2} - 1}{q-1} & 
       \displaystyle \frac{q^{m/2 - 1/2}}{q^{1/2} + 1} \\
\\
U,\bn & \displaystyle \frac{q^{m-1/2}}{q^{1/2} + 1} & 
       \displaystyle \frac{q^{m/2 - 1/2}}{q^{1/2} + 1} \\
\\
O,\bs & \displaystyle \frac{q^{m-1} - 1}{q-1} & 
       \displaystyle q^{m/2 - 1} \\
\\
O,\bn\ (\ q \mbox{\rm{ even }})  & \displaystyle q^{m-1} & 
       \displaystyle q^{m/2 - 1} \\
\\
O,\bn\ (\ q \mbox{\rm{ odd }})  & \displaystyle \frac{1}{2} q^{m-1} & 
       \displaystyle \frac{1}{2} q^{m/2 - 1/2} \\
\end{array}
$$
In particular, $N > q^{n-2}$.
\end{fact}

\pf
When $r = 0$, this follows immediately from Table~\ref{point.table} for 
all cases except odd-dimensional orthogonal spaces in even characteristic, 
in which case $\pi(W) = \AAA(m)$.  
The general case follows since 
$\pi(W) = \displaystyle \frac{q^r - 1}{q-1} + q^r \pi (W/R)$ for 
singular points and
$\pi(W) = q^r \pi (W/R)$ for nonsingular points.
%
%
\eop

\begin{fact}
Assume that $V$ is an even-dimensional unitary space or orthogonal 
space of type $+$.  
For $k \leq n/2 -1$, set 
$F_k = \AAA(n-k) + \BBB(n)$.
Set $F_{n/2} = 2\left( \frac{q^{n/2} - 1}{ q - 1 }\right)$.
If $q = q_0^2$, set $F^\ast = \frac{q_0^{n} - 1}{q_0 - 1}$.
\begin{enumerate}
\item
If $x$ is linear and $v(x) = k$, $k \leq n/2$, then 
$F(x) \leq F_k$.
\item
If $x$ is linear and $v(x) \geq n/2$, then 
$F(x) \leq F_{n/2}$.
\item
If $x$ is semilinear then $F(x) \leq F^\ast$.  
\end{enumerate}
\end{fact}

\pf 
This is a straightforward consequence of the previous statement. 
\eop

\begin{fact}
\label{geneigenspacefact}
Suppose $x$ preserves a non-degenerate sesquilinear or bilinear form on 
$V$, $X_\lambda = \ker (X-\lambda I)^n$, and
$X_\mu = \ker (X-\mu I)^n$.  
If $\lambda \bar{\mu} \neq 1$ then $X_\mu \subseteq X_\lambda^\perp$.
\end{fact}

\pf 
Argue by induction on $k + l$ that if $k$ and $l$ are positive 
integers, $v \in \ker (X- \lambda I )^k$, and
$w \in \ker (X- \mu I )^l$ then $\langle v , w \rangle = 0$.   
\eop


\subsection{Unitary and Orthogonal Groups}

\bigskip
To complete the proof of Theorem~\ref{grassmann.exceptions} we assume 
that $V$ 
admits a nondegenerate orthogonal or unitary form, and that the
action of $G$ is on the points of type $t$ in $V$.

To estimate $f(x) - q^{-v(x)}$ we bound $F(x)$ from above and $N$ from 
below.  
For a subspace $\U$ of $V$, 
let $\zzeta(\U) = \zzeta_t(\U)$ be the number of points of type $t$ in $\U$.
It is apparent that $F(x) = \sum \pi(E_\lambda)$ where $\{E_\lambda\}$ 
is the collection of eigenspaces for $x$.

\begin{lemma}
\label{Lemma2}
If $x$ acts linearly on $V$, then either
\begin{enumerate}
\item $f(x) - q^{-v(x)} < 1/100$ or
\item
\label{halfdimension.exception}
$V$ has even dimension, the action is on singular points,
and some eigenspace for the action of $x$ on
$V$ is a totally singular subspace of dimension $\dim V/2$.
\end{enumerate}
\end{lemma}

\pf
We have $F(x) = \sum \pi(E_\mu)$ where the sum is over the eigenspaces 
$E_\mu$ for the action of (some pull-back of) $x$ in the group of linear 
transformations of  $V$.

\bigskip
Suppose $v = v(x) \leq n/2$, and  
let $\lambda$ be the principal eigenvalue.
Then $\dim E_\lambda = n - v$. 
Let $X_\lambda$ be the 
corresponding generalized eigenspace, that is 
$X_\lambda = \ker (x-\lambda I)^n$, and set $w = \codim_V X_\lambda$.
Then $w \leq v$.

If $X_\lambda$ is totally singular, then $\dim X_\lambda \leq n/2$, and 
it follows that $X_\lambda = E_\lambda$, so the second alternative holds.  
We may therefore suppose that $X_\lambda$ is not totally singular.  
It follows from {\bf \ref{geneigenspacefact}\/} 
that $\lambda \bar\lambda = 1$ and that $E_\mu \subseteq X_\lambda^\perp$ 
whenever $\mu \neq \lambda$.

This implies that $$F(x) \leq \pi(E_\lambda) + \pi( X_\lambda^\perp).$$

Setting $r = \dim \rad(E_\lambda)$, we have 
$r \leq \codim_{X_\lambda}(E_\lambda) = v - w$.  
So $\dim X_\lambda^\perp = w \leq v-r$.
By {\bf \ref{point.count.lemma}\/}, 
$F(x) \leq \AAA(n-v) + \BBB(n-v+r) + \AAA(v-r) + \BBB(v-r)$ 
because $X_\lambda^\perp$ is non-degenerate.  

Since $N \geq \AAA(n) - \BBB(n)$ and $\AAA(n-v) \leq q^{-v}\AAA(n)$, 
it follows that
%
%
%
$$
\begin{array}{ccl}
F(x) - q^{-v}N  & \leq &
\BBB(n-v+r) + \AAA(v-r) + \BBB(v-r) + q^{-v} \BBB(n)  \\
  & = & 
\BBB(n-v+r) + \AAA(v-r) + \BBB(v-r) + \BBB(n-2v).
\end{array}
$$ 
because $\BBB(n) = K q^{n/2}$ where $K$ 
is independent of $n$.

Set $D(x) = F(x) - q^{-v} N$.  
We claim that $D(x) < (q+1) \BBB(n-2) + 2$.

We have shown that $D(x) \leq \phi(v,r)$ 
where $\phi(v,r) = \BBB(n-v+r) + \AAA(v-r) + \BBB(v-r) + \BBB(n-2v)$.
By elementary calculus, $\phi$ attains its maximum on the 
region $\{ (v,r) \ : \ 1 \leq v \leq n/2, 0 \leq r \leq v\}$ 
at $(1,1)$.
Therefore $D(x) \leq \phi(1,1) = \BBB(n) + \BBB(n-2) + \AAA(0) + \BBB(0)
< (q+1) \BBB(n-2) + 2$ since $\AAA(0) < 1$ and $\BBB(0) < 1$.

\newcommand\DDD{D}
Set $\DDD  = (q+1) \BBB(n-2) + 2$.  
It suffices to show that if $N \geq 10000$ then $\DDD/N < 1/100$.

In all cases, $N > q^{n-2}$ by {\bf \ref{point.count.lemma}\/}. 
When $V$ is unitary, 
$\BBB(n-2) = q^{(n-3)/2}/(q^{1/2}+1)$, and 
it is easy to see that $\DDD < q^{(n-2)/2}$. 
So $\DDD^2 < N$, which implies that $\DDD/N < 1/100$. 

We may therefore assume that $V$ is orthogonal.  
If either the action is on singular points or $q$ is even,
then 
$\BBB(n-2) = q^{n/2 -2}$,
and $\DDD < 
a^{n/2-1} (1 + q^{-1} + 2/(q^{n/2-1}) ) < 
8q^{(n-2)/2}/5$.
Therefore, $\DDD/N < \frac85 q^{- (n-2)/2}$, 
and $q^{(n-2)/2} < 160$.
This implies that $n \leq 16$ and $q \leq 11$ 
since $n \geq 6$.
Among the pairs $(n,q)$ of such values, the only ones 
for which both $q^{(n-2)/2} < 160$ and $N > 10000$ are
$(6,11), (7, 7) , (8,5) , ( 11, 3)$, and $(16,2)$.


In the nonsingular case when $q$ is odd,
$\BBB(n-2) = \frac12 q^{(n-3)/2}$,
so $\DDD \leq q^{(n-1)/2}(q+1)/2q + 2 < q^{(n-1)/2}$.
In this case, $N \geq \frac12 q^{n-1}(1 - q^{-(n-1)/2}) > \frac49 q^{(n-1)}$.
Therefore, $\DDD/N < \frac94 q^{- (n-1)/2}$, 
and $q^{(n-1)/2} < 225$.
This implies that $n \leq 10$ and $q \leq 7$.  
By inspection, the only pairs $( n , q )$ for which 
both $q^{(n-2)/2} < 160$ and $N > 10000$ are
$(6,8), (8,4) , ( 16, 2 )$.
A straightforward calculation shows that $f(x) - q^{-v} < 1/100$ in 
these cases.  


A straightforward calculation shows that $f(x) - q^{-v} < 1/100$ in 
these cases.  This shows that the result holds when
$v \leq n/2$.

\bigskip 
If $v \geq n/2 + 1$, 
then every eigenspace has dimension at most $n/2 -1$, and 
there are at most $q-1$ eigenspaces.
Therefore $F(x) \leq (q-1) \left( \displaystyle \frac{q^{n/2-1} -1}{q-1}\right)
< q^{n/2 -1} < \sqrt{N}$.  Since $N \geq 10000$ 
this implies that $F(x)/N < 1/100$.

\medskip
This leaves the case $v = (n+1)/2$, where $n$ is necessarily odd.  
Every eigenspace has dimension at most $(n-1)/2$,
so $F(x) \leq 2(q^{(n-1)/2} - 1)/(q-1) + 1 \leq F$ where $F = q^{(n-1)/2}$.
A short calculation, described below, shows that the conclusion holds in this case.

Suppose $V$ is unitary and 
set $q_0 = q^{1/2}$.  Then $N = 
\displaystyle \frac{q^{n-1/2} - 1}{q-1} - \frac{q^{(n-1)/2}}{q_0 + 1}$
in the singular case, and  
$N = \displaystyle \frac{q^{n-1/2} + q^{(n-1)/2}}{q_0 + 1}$ in 
the nonsingular case.
By computation, $F/N < 1/100$ when $(n,q_0) = (3,5), (5,3)$, or $(9,2)$.
Since $F/N$ is a decreasing function of both $q$ and $n$,
it follows that $n = 3, 5$, or $7$.  
Furthermore, $q_0 \leq 4$ when $n = 3$ and $q_0 = 2$ when $n = 5$ or $7$.
By inspection, $N < 10000$ in these cases.

In the orthogonal case, $q$ is necessarily odd because $n$ is odd.
We have
$N = 
\displaystyle \frac{q^{n-1} - 1}2$
in the singular case, and
$N = \displaystyle \frac{q^{n-1} \pm q^{(n-1)/2}}{2}$ in 
the nonsingular case.
By computation, $F/N < 1/100$ when $(n,q) = (7,7)$ or $(11,3)$.
Since $F/N$ is a decreasing function of both $q$ and $n$,
it follows that $n = 7$, or $9$.  
Furthermore, $q \leq 5$ when $n = 7$ and $q = 3$ when $n = 9$.
By inspection, $N < 10000$ in these cases.  

This completes the proof of Lemma~\ref{Lemma2}.  
\eop

\begin{lemma}
\label{Lemma3}
If $x$ acts semilinearly on $V$ then either
\begin{enumerate}
\item 
$f(x) - q^{-v(x)} < 1/100$ or
\item
\label{field.exception}
The dimension $n$ of $V$ is even, and $x$ has a totally singular 
eigenspace of dimension $n$ over $F_{q^{1/2}}$.
\end{enumerate}
\end{lemma}

\pf 
Assume that $x$ acts semilinearly on $V$ with $f(x) \geq 1/100$.
Let $d$ be the order of $x$ mod $PGL(V)$.
We claim that $d = 2$.

Suppose $d > 2$ and set $q_1 = q^{1/d}$.
Then the points fixed by $x$ must lie in an $n$-dimensional space 
over $GF(q_1)$, so
$F(x) \leq \displaystyle \psi(d) = \frac{q^{n/d} - 1}{q^{1/d} - 1}$.

If $V$ is orthogonal, and the action is on singular points,
then $\psi(d) < 1/100$ when $(d,q_1,n) = (3,2,8), (4,2,6), (3,3,6)$, or $(3,3,7)$.
If the action is on nonsingular points, then
Then $\psi(d) < 1/100$ when $(d,q_1,n) = (3,2,6), (3,3,6)$, or $(3,3,7)$.
For a given parity of $n$ and a given parity of $q$ the ratio
$\psi(d)$ is a decreasing function of $d,n$, and $q$.  
We have $N < 10000$ when $V$ is an orthogonal space of dimension 
$6$ over $F_8$, so $d=2$ in the orthogonal case.

In the unitary case, $q$ is necessarily a square.
We have $F(x)/N < 1/100$ when $(d,q_1,n) = (3,2^2,3), (3,2^2,4),
(4,2,3)$, or $(4,2,4)$, and the claim holds for the unitary case as well.

\smallskip
We have $d = 2$.
Let $E$ be the primary eigenspace for $x$. Then $E$ is an $\F_{q_0}$ space 
of dimension at most $n$ where $q_0 = q^{1/2}$.  
If $\dim E \leq n-1$, then 
$F(x) < (q_0^{n-1}-1)/(q_0-1) + 1$, and a short computation 
shows that $F(x)/N < 1/100$ whenever $N > 10000$.
Therefore $E$ has dimension $n$, and $F(x) \leq \psi(2)$.

We may assume that $\psi(2) - q_0^{-n}N  \geq N/100$, and in particular,
that $\psi(2) > N/100$.

Then $F(x) \leq \displaystyle \frac{q_0^n-1}{q_0 - 1}$, and 
the conditions $F(x)/N \geq 1/100$, $N \geq 10000$ imply that
$(n,q_0)$ is on one of the following lists.
%

Unitary groups, nonsingular action: 
%
$(8,2), (9,2), (4,5)$; \\
singular action:
$(8,2), (9,2), (6,3), (5,4), (4,7), (4,8), (4,9)$. \\
Orthogonal groups,
%
nonsingular action:
$(8,2), (6,3)$; \\
singular action:
$(10,2), (7,3), (6,4)$.

For $U_9(2^2)$, we have $\psi(2) \leq N/100 + q^{-n/2} N $.

For all other cases, the upper bounds in {\bf \ref{point.count.lemma}\/} 
imply that $F(x) < N/100 - q^{-n/2} N$ whenever $E$ is not totally singular.  
\eop

\begin{lemma}
\label{Lemma4}
If $f(x) - q^{-v} \geq 1/100$ for some $x \in G^\sharp$ 
then the action is on singular points and one of the following is true.
\begin{enumerate}
\item 
\label{u.alternative}
$V$ is unitary and $(n,q_0) \in \{ (4 , 7 ) , ( 4, 8 ) , ( 4, 9 ),
( 6 , 3 ) , ( 8 , 2 ) \}$.
\item 
\label{o.alternative}
$V$ is orthogonal of $+$ type and 
$(n,q) \in \{ (6,11) , ( 6, 13 ) , ( 6, 16 ),
( 8 , 5 ) , ( 10 , 4 ) \}$.
\end{enumerate}
\end{lemma}

\pf
The two previous lemmas show that it suffices to assume 
the action is on the singular points of an even-dimensional 
space $V$ and $V$ contains a totally singular subspace of dimension $\dim V/2$.

Suppose first that $V$ is a unitary space of dimension $2m$ over $\F_{q}$ where
$m \geq 2$ and $q = q_0^2$.  
Then $N = \displaystyle \frac{ (q_0^{2m} - 1 ) (q_0 ^ {2m - 1} + 1)}{q_0^2 - 1}$.

If $x$ has a totally singular eigenspace of dimension $m$, then the 
fixed points of $x$ are contained in the union of two subspaces of $V$ 
each of dimension $m$, so 
$F(x) \leq \displaystyle 2 \left( \frac{q_0^{2m} - 1 }{ q_0^2 - 1} \right)$.  
Otherwise,  $x$ has an eigenspace of dimension $2m$ over $\F_{q_0}$, 
and $F(x) \leq
\displaystyle \frac{q_0^{2m} - 1}{q_0 - 1}$.

In either case, 
$f(x) \leq \displaystyle \frac{{q_0}+1}{q_0^{2m-1} + 1 }$.  
By assumption, $f(x) \geq 1/100$.  
Therefore $q_0^{2m-2} < 100\left(\displaystyle \frac{{q_0}+1}{{q_0}} 
\right)\leq 150$.
Since $2^8 > 150$, it follows that $2m -2 < 8$.
Therefore $m \leq 4$. 
By inspection, one of the following holds:
$m = 2$, ${q_0} \leq 9$;
$m = 3$, ${q_0} \leq 3$; or
$m = 4$, ${q_0} = 2$.

By further inspection, $N < 10^4$ when $m = 2$, ${q_0}\leq 5$ and when 
$m = 3$, ${q_0} = 2$.
One of the conditions in \ref{u.alternative} must therefore hold.

\medskip
Now suppose $V$ is an orthogonal $+$ space of dimension $2m$ over $\F_q$, $m \geq 3$.
Then $N = \displaystyle \frac{ ( q^m - 1) ( q^{m-1} + 1 )}{ q - 1 }$.

Suppose $q = 2$. 
Then $x$ has a single eigenspace and $F(x) \leq q^m - 1$, so 
$f(x) \leq \displaystyle \frac{1}{q^{m-1} + 1}$.  Since $N > 10000$, we have $m \geq 8$.
Therefore $f(x) < 1/100$.  We may therefore assume that $q \geq 3$.

Suppose $x$ fixes a totally singular subspace of dimension $m$. 
Then the 
fixed points of $x$ are contained in the union of two subspaces of $V$ 
each of dimension $m$, so 
$F(x) \leq \displaystyle 2 \left( \frac{q^{m} - 1 }{ q - 1} \right)$.

We have $f(x) \leq \displaystyle \frac{2}{q^{m-1} + 1}$.
The assumption that $f(x) > 1/100$ implies that $q^{m-1} < 200$.
Since $q \geq 3$ and $3^5 > 200$, it follows that $m \leq 5$ and that 
one of the following holds:
$m = 3$, $q \leq 13$;
$m = 4$, $q \leq 5$; or
$m = 5$, $q = 3$.

By inspection, $N < 10^4$ when $m = 3$, $q \leq 9$, when $m = 4$, $q \leq 4$, and
when $m = 5$, $q = 3$, so one of the following must hold:
$2m = 6$ and $q = 11$ or $13$;  $2m = 8$ and $q = 5$.

\smallskip

If $x$ fixes a subspace of dimension $2m$ over $\F_{q^{1/2}}$, then $F(x) \leq
\displaystyle \frac{q^{m} - 1}{q^{1/2} - 1}$.
Therefore 
$f(x) \leq \displaystyle \frac{q^{1/2} + 1}{q^{m-1} + 1} $.
%

The condition $f(x) \geq 1/100$ implies that  $q^{m-1} + 1 \leq 100 (q_0 + 1)$, 
so $q_0^{2m-2} \leq 100(q_0 + 1) \leq 150 q_0$, and 
$q_0^{2m-3} \leq 150$.

We have $m \leq 5$ because $2^8 > 150$, and one of the following holds:
$m = 3$, $q_0 \leq 5$;
$m = 4$, $q_0 = 2$; or
$m = 5$, $q_0 = 2$.

By inspection, $N < 10^4$ for $m = 3$, $q_0 \leq 3$ and for $m = 4$, 
$q_0 = 2$.
Since $\frac{5+1}{25^2 + 1} < 1/100$, the case $m = 3$, $q_0 = 5$ does 
not satisfy the hypotheses.  This leaves the cases $2m = 6$, $q= 4^2$ 
and $2m = 10$, $q = 2^2$.
\eop 

\begin{lemma}
\label{newLemma5}
If one of the conclusions of Lemma~\ref{Lemma4} holds then 
$g(\ux) > 2$ whenever $\ux$ is a normalized generating tuple for $G$.
\end{lemma}

\pf
We consider the cases in turn.  
We assume throughout that $N > 10^4$, that $\ux$ is 
a normalized generating tuple for $G$, with signature $\ud$, and 
that $g(\ux) \leq 2$.  


By Theorem~\ref{basic.result}, it suffices to assume that 
there is an element $y$ involved in $\ux$ which violates 
Grassmann Condition $1/100$. 

\begin{stp}
\label{two.cases}
For some $i$, $\langle x_i \rangle$ 
contains an element $y$ such that one of the following is true.  
\begin{enumerate}
\item $y$ fixes two totally singular subspaces of dimension $n/2$.
\item $y$ is a semilinear map on $V$, $y$ has order $2$, and 
$y$ fixes a subspace of dimension $n$ over $F_{q^{1/2}}$.
\end{enumerate}
\end{stp}

\begin{stp} 
\label{ruleout237}
$\ux$ does not have signature $(2,3,7)$.
\end{stp}

\pf
If $\ux$ has signature $(2,3,7)$, then every element of $G$ 
must act linearly on $V$.  
By Step~\ref{two.cases}, $\ux$ must involve an element $y$ 
which has two totally singular eigenspaces of dimension $n/2$.  
No such element can violate Grassmann Condition $1/100$ when 
$G$ is of type $U_4(q_0^2)$, $U_6(3^2)$, $O_6(16)$, or $O_{10}(4)$.  
When $G$ is of type $U_8(2^2)$ no element can have two distinct 
totally singular eigenspaces.
In all other cases, the element of order $7$ can have at most 
one eigenvalue.  When $q = 11$ or $5$, the element of order $3$ 
can have at most one eigenvalue as well.  A short computation 
using fixed point estimates shows that $g(\ux) > 2$ in all cases.
\eop

\begin{stp}
$G$ is not of type $U_4(q^2)$.
\end{stp}

\pf 
$N = (q^2 + 1 ) ( q^3 + 1)$, and, 
for all $x \in G^\sharp$, 
we have $F(x) \leq F$ where $F = (q+1)(q^2 + 1)$. 
The Riemann-Hurwitz Formula implies that $A(\ud) \leq (2N + 2 )/(N - F)$.
Since $\ud \neq (2,3,7)$, we must have $q=7$, $\ud = (2,3,8)$.
In this case, $v(x_3^2) > 1$, and $x_3^2$ must act linearly on $V$, 
so $F(x_3) \leq F(x_3^2) \leq 2(q^2 + 1)$.  
By computation, $g(\ux) > 2$.
\eop


\begin{stp}
\label{O6step}
$G$ is not of type $O_6^+(q)$.
\end{stp}

\pf
We have $N = (q^3-1)(q^2+1)/(q-1)$.
Set $F_1 = (q^4-1)/(q-1) + q^2$,
$F_2 = (q^3-1)/(q-1) + q^2$, and
$F_3 = 2(q^3-1)/(q-1)$.
When $q = 16$, set $q_0 = 4$ 
and $F^\ast = (q^3-1)/(q_0-1)$.

Then $F(x) \leq F_1$ for all $x \in G^\sharp$, 
and the Riemann-Hurwitz Formula implies that 
$A(\ud) \leq (2N + 2)/(N- F_1)$.
It follows that $\ud$ is one of the following:
$(2,3,d)$;
$(2,4,d)$, $d \leq 29$;
$(2,5,d)$, $d \leq 11$;
$(2,6,d)$, $d \leq 8$;
$(2,7,7)$;
$(3,3,d)$, $d \leq 8$;
$(3,4,4)$; or 
$(2,2,2,3)$.

By inspection, if $\cB$ is the set of all elements in $G^\sharp$ for 
which $F(x) > \max(F_2,F^\ast)$ then 
$\sum \displaystyle \frac{|\{ \cB \cap \langle x_i \rangle \}|}{d_i} < 1$.
Set $F' = \max(F_2, F_3, F^\ast)$.

The Riemann-Hurwitz Formula now implies that 
$A(\ud) \leq (2N + 2 + 1(F_1-F'))/(N- F')$, whence
$\ud$ is one of the following:
$(2,3,d)$, $d \leq 19$;
$(2,4,d)$, $d \leq 7$;
$(2,5,5)$; or
$(3,4,4)$.

Inspecting this list it follows that
$\sum \displaystyle \frac{|\{ \cB \cap \langle x_i \rangle \}|}{d_i} \leq 1/2$, 
so $A(\ud) \leq (2N + 2 + \frac12(F_1-F'))/(N- F')$, which further 
reduces the possible signatures. 
Further iterations of this procedure show that $\ux$ must have 
signature $(2,3,7)$, which was already ruled out by Step~\ref{ruleout237}.
%
\eop

\begin{stp}
\label{U6step}
$G$ is not of type $U_6(3^2)$.
\end{stp}

\pf 
In this case, using 
$N = (q_0^6-1)(q_0^5 + 1) / ( q - 1 )$,
$F_1 = (q_0^9 -1)/ ( q - 1 ) + q_0^5/(q_0 + 1)$,
$F_2 = (q_0^7 -1)/ ( q - 1 ) + q_0^5/(q_0 + 1)$,
$F_3 = 2(q^3 - 1 )/( q - 1 )$, 
and $F^\ast = (q_0^6 - 1) / ( q_0 - 1)$, 
a short modification of the analysis in the previous step 
again reduces to the case $\ud = (2,3,7)$, which was treated earlier. 
\eop

\begin{stp}
$G$ is not of type $O_8^+(5)$.
\end{stp}

\pf 
The argument of the previous two steps shows that either 
$\ud = (2, 3, d)$ for some $d$  or
$\ud \in \{ (2, 4, \leq 8), (2, 5, \leq 6), ( 3, 3, \leq 5) \}$.

In the former situation, the contribution of elements 
having $v(y) = 1$ is less than $2/3$, and it follows that $d < 200$,
whence the contribution is less than $5/12$.  Continuing in 
this way shows that no tuple $\ux$ can have $g(\ux) \leq 2$.

In the remaining cases, bounding the contributions from elements 
with $v(y) \leq 2$ leads to the same conclusion.
\eop

\begin{stp}
\label{U8step}
$G$ is not of type $U_8(2^2)$.
\end{stp}

\pf
The argument of the previous steps shows that either 
$\ud = (2, 3, d)$ or $(2, 4, d)$ for some $d$  or
$\ud \in \{ (2, 5, \leq 17), (2, 6, \leq 10), (2, 7, \leq 8), 
( 3, 3, \leq 10) , ( 3, 4, \leq 5) , ( 2, 2, 2, 3 ) \}$.

If $\beta_1$ is the contribution to $A(\ud)$ from elements 
having $v(y) = 1$ and 
$\beta_2$ is the contribution from elements with $v(y) \leq 2$,
then the Riemann-Hurwitz Formula implies that 
$A(\ud) \leq  (2N + 2 - \beta_1(F - F') - \beta_2(F' - F'') )
/ (N - F'')$, 
where $F(x) \leq F$ for all $x \in G^\sharp$,
$F(x) \leq F'$ for all $x$ with $v(x) > 1$, and
$F(x) \leq F''$ for all $x$ with $v(x) > 2$.  

Using this criterion eliminates the 
individual cases other than $(2,3,d)$, $(2,4,d)$.
Using estimates for indexes, this reduces to $(2,3, \leq 19)$,
or $(2,4, \leq 9)$.

In the $(2,3,d)$ case, we have $\Ind(x_2) \geq \frac23(N - 
2 (q_0^4-1)(q_0^3+1)/(q-1) )$ because $v(x_2) \geq 4$ and 
the eigenspaces for $x_2$, an element of order $3$, must 
be nondegenerate.  Bounding $v(x_3^k)$, and hence $F(x_3^k)$, 
for $k = 1,2,3,4$, shows that $g(\ux) > 2 $ for for all choices 
of $d$.

In the $(2,4,d)$ case, consideration of the subcases $v(x_2^2) = 1$, 
$v(x_2^2) > 1$ leads to the same conclusion.
\eop

\begin{stp}
$G$ is not of type $O_{10}^+(4)$.
\end{stp}

\pf
The argument in Step~\ref{O6step} 
shows that either 
$\ud = (2, 3, d)$ or $(2, 4, d)$ for some $d$  or
$\ud \in \{ (2, 5, \leq 12), (2, 6, \leq 9), (2, 7, 7), 
( 3, 3, \leq 9) , ( 3, 4, \leq 5) , ( 2, 2, 2, 3 ) \}$.

Its extension in 
Step~\ref{U8step} reduces to the earlier treated case $\ud = (2,3,7)$.
%
%
%
%
\eop 

Combining Lemmas~\ref{Lemma4} and 
\ref{newLemma5} we have the following result.

\begin{prop}
\label{OUprop}
If $G$ is unitary or orthogonal and $G$ violates Grassmann 
Condition $1/100$ then $g(\ux) > 2$ for every normalized generating 
tuple for $G$.
\end{prop}

Propositions~\ref{linear}, \ref{Sphypprop}, and \ref{OUprop}
establish Theorem~\ref{grassmann.exceptions}.

\section{Proof of Theorem~\ref{touch.up}}

We assume here that $\ux$ and $V$ satisfy one of the conditions
listed in Table~\ref{pnd.table}.
Suppose $\Omega$ is a primitive $G$-set of [projective] points in $V$
with $\order{\Omega} \geq 10000$.

That is, one of the following is true where $n_p = \dim_{\F_p}(V)$.
      \begin{enumerate}
         \item $\ux$ has signature $(2,3,7)$ and one of the following holds.
            \begin{enumerate}
               \item $p = 11$ and $n_p = 5$ or $6$.
               \item $p = 7$ and $n_p = 6$.
               \item $p = 5$ and $n_p = 7,8$, or $9$.
               \item $p = 3$ and $n_p = 12$.
               \item $p = 2$ and $14 \leq n_p \leq 21$.
            \end{enumerate}
         \item $\ux$ has signature $(2,3,8)$, $p = 3$, and $n_p = 10$.
         \item $\ux$ has signature $(2,4,5)$, $p = 2$, and $n_p = 16$
     Furthermore $v_p(x_1) = 4$, $v_p(x_2) = 12$, 
                 and $v_p(x_3) = 16$.
      \end{enumerate}

Then $V$ is an $n_q$-dimensional $\F_q$-module where $q^{n_q} = p^n$.
and $n_q$ and $q$ satisfy the conditions listed for point actions.

\begin{fact}
The number $CP(X,n,q,t)$ of $t$-points in a classical $n$-space of type X over $GF(q)$
is given in Table~\ref{point.table}.
\end{fact}
 
\pf
See \cite{FM1}.
\eop

We calculate a lower bound for $g(\ux)$ in each of the cases using the following
lemma.

\begin{lemma}
\begin{enumerate}
\item If $\ud = (2,3,7)$, then $v(x_1) \geq n/3$,
$v(x_2) \geq n/2$, and $v(x_3) \geq n/2$.
\item If $\ud = (2,3,8)$ then $v(x_1) \geq n/3$,
$v(x_3) \geq n/2$, $v(x_2) \geq n/2$, $v(x_2^2) \geq n/2$, and $v(x_2^4) \geq n/5$. 
\item If $\ud = (2,4,5)$, then $v(x_1) \geq n/4$,
$v(x_2) \geq n/2$, $v(x_2^2) \geq n/4$, and $v(x_3) \geq n/2$.
\item 
The number of $t$-points in an $n$-space with radical of dimension $r$
of type $X$ over $\F_q$ is $(q^r-1)/(q-1) + q^r CP(X,n-r,q,t)$ for singular points and
$q^r CP(X,n-r,q,t)$ for non-singular points.
\item 
\label{fifth.statement}
Assume that $q$ is even.
Let $G = O(2m+1,q) \cong Sp(2m,q)$ act on the $2m+1$-dimensional orthogonal space $V$,
where $V$ has a $1$-dimensional radical $R$.
If $x$ is a linear transformation in $G$ then $x$ fixes at most
$q^m(q^{m-v(x)}+1)/2$ complements to $R$ of each type.
\item If $W$ is a space of codimension $v$ in the non-degenerate space $V$
then $\dim \rad W \leq v$.
\item 
\label{seventh.statement}
Let $\Fix_2(x)$ be the number of fixed points of $x$ lying outside its principal
eigenspace.  Set $v = v(x)$.
Then 
\begin{enumerate} 
\item 
\label{seventh.1}
If $(o(x),q-1) = 1$ then $\Fix_2(x)=0$.
\item 
\label{seventh.2}
$\Fix_2(x) = 0$ in case of type $S$.
\item 
\label{seventh.3}
If $2v \leq n$ then $\Fix_2(x)$ is bounded by the number of type $t$ points
in some $v$-dimensional space.
\item 
\label{seventh.4}
If $(o(x),q-1) = d_0$ and every $n-v$-dimensional space contains at most
$M$ points then $\Fix_2(x) \leq (d_0-1)M$.
\end{enumerate}
\item 
\label{eighth}
If $\Fix(x^j) \leq F_j$ for all positive powers of $x$, then 
$$\Ind(x) \geq \frac{d-1}{d}N - \frac1d \left( \sum_{k|d,k<d} \phi(\frac{d}{k})F_k \right)$$
\item 
\label{ninth}
If $\Ind x_i \geq H_i$ for all $i$ then $g(\ux) \geq \frac12\sum H_i -N + 1$.
\end{enumerate}
\end{lemma}

\pf The first three statements follow from Lemma~\ref{translation.lemma}.

The fourth statement is a straightforward count of points in $R \oplus W$ where
$R$ is totally singular of dimension $r$ and $W$ is non-degenerate.

Statement \ref{fifth.statement} follows from a straightforward 
calculation, as in the proof of Proposition 8.1 of \cite{FM1}.

The next statement is clear because $ \rad W \subseteq W^\perp$.

To prove \ref{seventh.statement}, note that
the principal eigenspace of $x$ has dimension $n-v$, and
every fixed point of $x$ lying outside the principal eigenspace must lie in
an eigenspace of dimension at most $n-v$.  

All eigenvalues of $x$ must
have order dividing both $o(x)$ and $q-1$, so there are at most $d_0 = (o(x),q-1)$
eigenvalues in toto.  Statements \ref{seventh.1} and \ref{seventh.4}
now follow immediately.

In type $S$ only the eigenvalue $\lambda = 1$ corresponds to fixed points,
so statement \ref{seventh.2} holds.

The total dimension of all secondary eigenspaces is at most $v$, and all
secondary fixed points of $x$ lie in the direct sum of such subspaces.
Statement \ref{seventh.3} follows.

Statements \ref{eighth} and \ref{ninth} follow easily from the Cauchy-Frobenius
and Riemann-Hurwitz Formulas, respectively.
\eop

\begin{center}
\begin{table}
\caption{Number of $t$-points in classical $n$-space of type $X$ over $\F_q$.}
\begin{tabular}{cccc}
\label{point.table}
$X$ & condition & $t$ & $CP(X,n,q,t)$ \\
\\
$L$ & & & $\displaystyle \frac{q^n-1}{q-1}$ \\
\\
$O^\eps$ & $n = 2m$ & singular  & $\displaystyle\frac{(q^m-\eps1)(q^{m-1}+\eps1)}{q-1}$ \\
\\
$O^\eps$ & $n = 2m$ & $\delta$  & $\displaystyle\frac{(2,q)}2 {(q^m-\eps1)q^{m-1}} $ \\
\\
$O^\eps$ & $n = 2m+1$ & singular  & $\displaystyle\frac{q^{2m}-1}{q-1}$ \\
\\
$O^\eps$ & $n = 2m+1$ & $\delta$  & $\displaystyle\frac{q^m(q^{m}-\eps\delta)}{2}$ \\
\\
$U$ & $q = q_0^2$ & singular & $\displaystyle \frac{(q_0^{n} - (-1)^{n})(q_0^{n-1} + (-1)^{n})}{q-1}$ \\
\\
$U$ & $q = q_0^2$ & non-singular & $\displaystyle \frac{(q_0^{n}-(-1)^{n})q_0^{n-1}}{q_0 + 1}$ \\
\\
$S$ & $n = 2m,q$ even & $\eps$ hyperplane & $\displaystyle \frac{q^{m}(q^{m}+\eps1)}{2}$ \\
\end{tabular}
\end{table}
\end{center}

In all cases except $L_{14}(2)$ acting on the points in its natural module
and $U_8(2^2)$ acting on singular points the lower bound is larger than 2.

However, in those cases, we use the following additional facts:
\begin{enumerate}
\item If $x$ has order $7$
and acts as a linear transformation over $\F_2$ or $\F_4$ then $x$ has a single
eigenspace and $3|v(x)$.
\item If $x$ has order $3$
and acts as a linear transformation over $\F_2$ then $x$ has a single
eigenspace and $2|v(x)$.
\end{enumerate}

Using these additional facts, it is easy to establish the 
following lemma and complete the proof of Theorem~\ref{touch.up}.
\begin{lemma}
If $\ud = (2,3,7)$ and the action is either $L_{14}(2)$ on points or
$U_{8}(2^2)$ on singular points, then the genus is at least $20$.
\end{lemma}

\pf
Suppose $G = L_{14}(2)$. 
Then $x_i$ has only one eigenspace for $i = 1, 2, 3$,
$2 | v(x_2)$, and $3 | v(x_3)$.  It follows that
$v_1 \geq 5$, $v_2 \geq 8$, and $v_3 \geq 9$. 
Furthermore, 
$\Ind(x_1) \geq \frac12(2^{14}-2^9) = 7936$,
$\Ind(x_2) \geq \frac23(2^{14}-2^6) = 10880$, and
$\Ind(x_3) \geq \frac67(2^{14}-2^5) = 14016$.
This implies that $g(\ux) > 30$.

Suppose $G = U_8(2^2)$.
Then $x_1$ and $x_3$ have at most one eigenspace,
and $3 | v(x_3)$.  
We have $v_1 \geq 3$, $v_2 \geq 4$, and $v_3 = 6$,
and it follows that $g(\ux) > 2$.
\eop

\bibliographystyle{alpha}
\bibliography{publist}

\end{document}